\documentclass[a4 paper,12pt]{article}
\usepackage[T1]{fontenc}
\usepackage[latin1]{inputenc}
\usepackage[german,english]{babel}
\usepackage{graphicx}
\usepackage{mathrsfs}
\usepackage{amsfonts}
\usepackage{amsxtra}
\usepackage{amssymb}
\usepackage{eufrak}
\usepackage{pst-all}
\usepackage{multido}
\usepackage{amsmath}
\usepackage{color}
\usepackage{makeidx}
\usepackage{multicol}
\usepackage[absolut]{overpic}

\usepackage[rm,bf,tiny,center]{titlesec}
\titlelabel{\thetitle.\enspace}
\newcommand{\ger}[1]{\mathfrak{#1}}

\newcommand{\OO}{\mathcal{O}}

\newcommand{\A}{\mathcal{A}}
\newcommand{\ma}{\leqslant}
\newcommand{\De}{\Delta}
\newcommand{\CC}{\mathbb{C}}

\newcommand{\NN}{\mathbb{N}}

\newcommand{\II}{\mathcal{I}}

\DeclareMathOperator{\add}{add}
\DeclareMathOperator{\dimm}{dim}
\DeclareMathOperator{\Hom}{Hom}
\DeclareMathOperator{\Kern}{ker}

\DeclareMathOperator{\End}{End}
\DeclareMathOperator{\GL}{GL}
\DeclareMathOperator{\topp}{top}
\DeclareMathOperator{\rad}{rad}
\DeclareMathOperator{\modd}{mod}
\DeclareMathOperator{\soc}{soc}

\DeclareMathOperator{\Ext}{Ext}
\DeclareMathOperator{\im}{im}
\DeclareMathOperator{\spann}{span}

\newtheorem{num}{}[section]

\usepackage{geometry}

\usepackage{fancyhdr}
\begin{document}





\begin{center}
\begin{large}\textbf{1-quasi-hereditary algebras:\\ Examples and invariant submodules of projectives} \end{large}
\end{center}
\begin{center}
Daiva Pu\v{c}inskait\.{e}
\end{center}
\begin{abstract}
In \cite{P} it was shown that every 1-quasi-hereditary algebra $A$ affords a particular basis which is related to the partial order $\ma$ on the set of simple $A$-modules.
In this paper we show that  the modules generated by these  basis-elements are also modules over the endomorphism algebra of some projective indecomposable modules.  
 In case the   Ringel-dual of 1-quasi-hereditary algebra is also 1-quasi-hereditary, all local $\De$-good submodules of projective indecomposable modules are also $\End_A(P)^{op}$-modules for the projective-injective indecomposable module $P$.
\end{abstract}

\begin{center}
\textbf{Introduction}
\end{center}
The class of (basic)  1-quasi-hereditary algebras, introduced in \cite{P}  is characterized by the fact that  
      all possible non-zero Jordan-Hölder multiplicities of standard modules as well as   $\De$-good multiplicities  of indecomposable projectives   are equal to 1.       
      Many   factor algebras (related to a  saturated subsets)  of an algebra $\A$ associated  to a block  of   category $\OO(\ger{g})$  of  a  semisimple  $\CC$-Lie algebra  $\ger{g}$ are 1-quasi-hereditary. 
      
   The class of    1-quasi-hereditary algebras 
 has a non-empty intersection with  some  other subclasses of quasi-hereditary algebras: BGG-algebras \cite{Chan}, quasi-hereditary algebras having Borel subalgebras \cite{K}, Ringel self-dual algebras etc.. 
However,  1-quasi-hereditary algebras are in general not BGG-algebras and the  class of 1-quasi-hereditary algebras  is    not closed under Ringel duality. The selected examples in this paper  serve to illustrate this.

Several properties of 1-quasi-hereditary algebras only  depend on the related partial order ($\Ext$-quiver,  good-filtrations, etc. see \cite{P}). In particular,  a  1-quasi-hereditary  algebra $(A,\ma)$ has a $K$-basis $\ger{B}(A)=\dot{\bigcup}_{j\in Q_0(A)}\ger{B}_j(A)$ containing distinguished  paths, which are  linked  only to $\ma$.  The paths in  $\ger{B}_j(A)$  form a $K$-basis of the projective  indecomposable $A$-module  $P_A(j)$ for any $i\in Q_0$ (see \cite[Theorem 3.2]{P}). This basis plays an important role for the structure of the  endomorphism algebras of the projective $A$-modules.
\\

\hspace*{-5mm}\textbf{Theorem A}. \textit{Let $(A, \ma)$ be a   1-quasi-hereditary algebra and $j\in Q_0(A)$.  The submodules of $P_A(j)$ generated by the paths in   $\ger{B}_j(A)$ are   $ \End_A(P_A(j))^{op}$-modules. 
}
\\

Any 1-quasi-hereditary algebra $A$ has (up to isomorphism) a unique projective-injective indecomposable module $P(A)$ (see \cite[Lemma 2.1]{P}). 
In  case the Ringel-dual $R(A)$ of  $A$ is also 1-quasi-hereditary, there exist  elements in $P(A)$  which generate  modules  satisfying    remarkable properties.
\\

\hspace*{-5mm}\textbf{Theorem B}. \textit{Let $(A, \ma)$ and $(R(A),\geqslant)$ be    1-quasi-hereditary algebras.  There exists a $K$-basis $\ger{B}$   of the projective-injective indecomposable $A$-module  $P(A)$,  such that
\begin{itemize}
	\item[(1)] The set $\left\{A\cdot \ger{b}\mid \ger{b}\in \ger{B}\right\}$ 
	is the set of  local,  $\De$-good submodules of all projective indecomposable $A$-modules.
\\[30pt]
 \hspace*{-0.8cm}\begin{minipage}{6cm}
 \hrule
 \vspace{0.5cm} \textcolor{black}{ }
\end{minipage}
\\[-4mm]
\hspace*{-0.7cm}\begin{footnotesize}\rm{Partly supported by the D.F.G. priority program SPP 1388 ``Darstellungstheorie''.}\end{footnotesize}
	\item[(2)] For every $\ger{b}\in \ger{B}$  the space $A\cdot \ger{b}$ is an $\End_A(P(A))^{op}$-submodule of $P(A)$.
\end{itemize}
}

The paper is  organised    in the following way: In  Section 1,  the definition of 1-quasi-hereditary algebra is recalled,  and we  present some examples of 1-quasi-hereditary algebras. 
In particular,  we  define a 1-quasi-hereditary algebras $A_n(C)$ for some $C\in \GL_{n-2}(K)$ for $n\geq 3$.  We show   that these algebras are not BGG-algebras in general and that the  Ringel-dual of $A_n(C)$  is 1-quasi-hereditary,  namely $R(A_n(C)) \cong A_n(C^{-1})$. 
Moreover,  we give  an example of a  1-quasi-hereditary algebra, whose Ringel-dual is not 1-quasi-hereditary. Using this algebra we 
illustrate the connection between the Jordan-Hölder-filtrations of standard resp. costandard and good-filtrations of injectives resp. projectives indecomposable  modules described in \cite[Sec. 4]{P}.

Section 2 and  3 is devoted to the proof of Theorems  A  and B respectively. 
These are based on order reversing and preserving bijections  between  certain subsets of $Q_0(A)$ with respect to $\ma$ and the  lattice of modules generated by $p\in \ger{B}(A)$.
The connections between the modules and the partial order $\ma$ given  in Theorem A and B will be illustrated  using the algebra associated to the regular block of category  $\OO(\ger{sl}_3(\CC))$.

\section{Preliminaries  and examples of 1-quasi-hereditary algebras}
\begin{small} Throughout the paper,  any algebra $\A$ is a  finite dimensional, basic  $K$-algebra over an algebraically closed field $K$ given   by a quiver $Q(\A)$ and relations $\II(\A)$.    The vertices of $Q(\A)$    are    parameterized   by the natural numbers. Any  $\A$-module  is a  finite dimensional  left module and for any $i\in Q_0(\A)$ we denote  by $P(i)$, $I(i)$ and $S(i)$  corresponding  projective, injective and simple $\A$-module,  respectively.   The product of arrows $(i\to j)$ and $(k\to i)$ is given by  $(k\to i\to j)=(i\to j )\cdot (k\to i)$. If $\left|\left\{i\stackrel{\alpha}{\rightarrow}j \mid \alpha \in Q_1(\A)\right\}\right|\leq 1$ for all $i,j\in Q_0(\A)$, then a path $p=(i_1\to i_2\to \cdots \to i_m)$  in $Q(\A)$  (and also  in $\A$) we denote  by $(i_1i_2\ldots i_m)$.  
For the corresponding  $\A$-map $f_p:P(i_m)\to P(i_1)$ given by   $f(e_{i_m})=p$ holds $f_p=f_{(i_1 i_2)}\circ f_{(i_2 i_3)}\circ\cdots\circ f_{(i_{m-1}i_m)}$. The $\A$-module generated by $p$ we denote by $\left\langle p\right)$, i.e.  $\left\langle p\right) = \A\cdot p =\im(f_p)$.
By 
$[M:S(i)]=\dimm_KM_i$ we denote the Jordan-Hölder multiplicity of $S(i)$ in  an $\A$-module $M$.\end{small}
\\[-2mm]

    The equivalent definition  of quasi-hereditary   algebras introduced by Cline-Parshall-Scott \cite{CPS} 
  is given by Dlab and Ringel in \cite{DRin}. The following is a brief review of  relevant terminology and notations: Let $(Q_0(\A),\ma)$ be  a  partially ordered set.   
For  every   $i\in Q_0$ the  \textit{standard}  module    $\De(i)$  is   the largest  factor module    of 
    $P(i)$  such that      $[\De(i):S(k)]=0$ for all $k\in Q_0$ with  $k\not\ma i$,  resp.   the  \textit{costandard}  module   $\nabla(i)$ is  the largest submodule of 
    $I(i)$ such that   $[\nabla(i):S(k)]=0$ for all $k\in Q_0$ with   $k\not\ma i$. 
   We denote   by $\ger{F}(\De)$   the full subcategory of $\modd \A$  consisting of the modules having a filtration  such that each subquotient is isomorphic to a standard module.    The modules   in  $\ger{F}(\De)$  are called    \textit{$\De$-good}  and  these   filtrations are   \textit{$\De$-good filtrations}      (resp.  \textit{$\nabla$-good} modules have     \textit{$\nabla$-good filtrations} and belongs to $\ger{F}(\nabla)$).   For  $M\in \ger{F}(\De)$  we   denote  by $\left(M:\De(i)\right)$  the (well-defined) number of  subquotients  isomorphic to  $\De(i)$ in some $\De$-good filtration  of $M$ (resp. $\nabla(i)$ appears  $(M:\nabla(i))$ times   in some  $\nabla$-good filtration of $M\in \ger{F}(\nabla)$).
\\[8pt]
\hspace*{2mm}\parbox{16cm}{
The algebra $\A=\left(KQ/\II, \ma\right)$ is   \textit{quasi-hereditary}    if  for all $i,k\in Q_0$ the following   holds: 
\begin{itemize}
	\item   $[\De(i):S(i)]=1$,
	\item $P(i)$ is a  $\De$-good module   with $\left(P(i):\De(k)\right)=0$ for all $k\not\geqslant i$ and $\left(P(i):\De(i)\right)=1$.
\end{itemize}
     }

\begin{num}\normalfont{\textbf{Definition.}}
A quasi-hereditary algebra $A=(KQ/\II, \leqslant)$
is called 
\textit{1-quasi-hereditary}
if for all $i,j\in Q_0=\left\{1, \ldots , n\right\}$ the
following conditions are satisfied:
\begin{itemize}
    \item[(1)] There is a smallest and a largest element with respect to  $\ma$, 
    \\ without loss of generality we will assume  them to be $1$ resp. $n$,
    \item[(2)] $[\Delta(i):S(j)]=\big(P(j):\Delta(i)\big)=1$ for $j\leqslant i$, 
\item[(3)] $\soc P(i) \cong \topp I(i) \cong S(1)$,
\item[(4)] $\Delta(i) \hookrightarrow \Delta(n)$ and $\nabla(n)\twoheadrightarrow
\nabla(i). $
\end{itemize} \label{def1qh}
\end{num}
The projective indecomposable module of a 1-quasi-hereditary algebra $A$ which corresponds to the minimal vertex $1$ is also injective with $P(1)\cong I(1)$ (see \cite[Lemma 2.1]{P}).
Axiom (3) shows that any projective indecomposable  $A$-module can be considered as a submodule of $P(1)$. Axiom (2) implies that for any $i\in Q_0(A)$ there exist  a   uniquely  determined submodule $M(i)$ of $P(1)$  with $M(i)\cong P(i)$. Consequently, for all $i,j\in Q_0$ with $j\ma i$ there exists  a   uniquely   determined submodule of $P(j)$ isomorphic to $P(i)$ (see  \cite[Lemma 2.2]{P}). We will often make use of this fact in the following. 
\\

\textbf{Example 1.} Let $\A$ be an algebra associated to a block of the  category    $\mathcal{O}(\ger{g})$  of    a complex semisimple  Lie algebra  $\ger{g}$  defined in  \cite{BGG}. Then $A$  is 1-quasi-hereditary if  $\text{rank}(\ger{g})\leq 2$. The quivers and relations of these   algebras  are to be found
 in \cite{St1}. To illustrate some statements   we use the algebra  corresponding  to 
a regular block of  category 
$\OO(\mathfrak{sl}_3(\CC))$ (see also \cite{Ma}):   
\\
\\

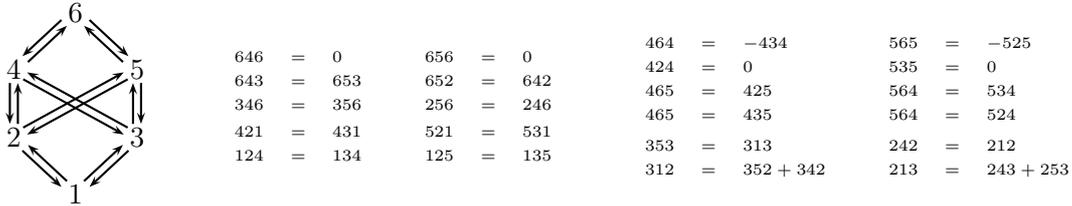
\psset{xunit=2.7mm,yunit=2.3mm,runit=4mm}
\begin{pspicture}(0,-10)(0,0)
\rput(5,1){\begin{small}
              \rput(0,0.5){\rnode{0}{6}}
              \rput(-3,-2.8){\rnode{1}{4}}
              \rput(3,-2.8){\rnode{2}{5}}
              \rput(-3,-6.7){\rnode{3}{2}}
              \rput(3,-6.7){\rnode{4}{3}}
              \rput(0,-10){\rnode{5}{1}}  \end{small}          }
\rput(33,-4){\begin{tiny}$
\begin{array}{ccl} 646& = &0\\[2pt]
 643& = &653\\[2pt]
 346&=&356\\[3pt]
421& = & 431\\[2pt]
124 & = & 134\\
\end{array}
$ \hspace{3mm} $
\begin{array}{ccl}
 656& = &0\\[2pt]
 652& = &642\\[2pt]
 256&=&246\\[3pt]
521& = & 531\\[2pt]
125 & = & 135\\ 
\end{array}$ \hspace{7mm} $
\begin{array}{ccl} 
 464& = & -434\\[2pt]
 424& = & 0\\[2pt]
 465& = & 425\\[2pt]
465 & = & 435\\[5pt]
353 & = & 313\\[2pt]
312 & = & 352+342\\
\end{array}
$ \hspace{3mm} $
\begin{array}{ccl}
 565& = & -525\\[2pt]
 535& = & 0\\[2pt]
 564& = & 534\\[2pt]
564 & = & 524\\[5pt]
242 & = & 212\\[2pt]
213 & = & 243+253\\
\end{array}$
 \end{tiny}}             
 \psset{nodesep=1.5pt,offset=1.5pt,arrows=<-}
              \ncline{0}{1}
              \ncline{1}{0}
              \ncline{1}{3}
               \ncline{3}{1}
               \ncline{1}{4}
               \ncline{4}{1}
               \ncline{3}{5}
               \ncline{5}{3}
              \psset{nodesep=1.5pt,offset=1.5pt,arrows=->}
              \ncline{2}{4}
              \ncline{4}{2}
              \ncline{0}{2}
              \ncline{2}{0}
              \ncline{2}{3}
              \ncline{3}{2}
              \ncline{4}{5}
              \ncline{5}{4}
              
\end{pspicture}

 The quiver and relations of the  algebra  $\A$ corresponding to a regular block of $\mathcal{O}(\ger{sl}_4(\CC))$  are  calculated  in   \cite{St1}  (in this  notations we have   $24\ma i\ma 1$  for all $i\in Q_0(\A)$). The algebra $\A$  is not 1-quasi-hereditary, since $[\De(3):S(16)]=2$. 	However,  the factor algebra $\A/\left(\A \epsilon(j)\A\right)$ with $\epsilon(j)=\sum_{i\not\ma j}e_i$  is 1-quasi-hereditary  for every $j\in \left\{i\in Q_0\mid i\ma x \text{ for some } x=5,6,7,9\right\}$. 
 \\

\textbf{Example 2.} Dlab, Heath and Marko described in  \cite{DHM} quasi-hereditary algebras which are obtained  in the following way: Let $B$ be  a commutative local self-injective $K$-algebra, $\dimm_KB=n$. Let $\ger{X}=\left\{\mathcal{X}(\lambda)\mid \lambda\in \Lambda\right\}$ be a set of local ideals of $B$  with  $B=\mathcal{X}(\lambda_1)\in \ger{X}$  indexed by a finite partially ordered set $\Lambda$ reflecting inclusions:  $\mathcal{X}(\lambda')\subset \mathcal{X}(\lambda)$ if and only if $\lambda'>\lambda$. Then $ A=\End_B\left(\bigoplus_{\lambda\in \Lambda}\mathcal{X}(\lambda)\right)$ is a quasi-hereditary algebra with respect to $(\Lambda,\leqslant)$ if and only if $\left|\Lambda\right|=n$ and $ \rad \mathcal{X}(\lambda)=\sum_{\lambda < \mu}\mathcal{X}(\mu)$ for every   $\mathcal{X}(\lambda)\in \ger{X}$.

Quasi-hereditary  algebras obtained in this way satisfy  all  axioms  for  1-quasi-hereditary algebra (see  \cite[Section 4]{DHM}). 
\\

The   algebras  in the examples 1 and 2   are    BGG-algebras: A quasi-hereditary algebra   $\A=(KQ/\II,\ma)$ with  a duality  functor  $\delta$ on $\modd \A$ [$\delta$ is a contravariant, exact, additive functor   such that $\delta \cdot \delta$ is the identity on $\modd A$ and $\delta$ induces a $K$-map on the  $K$-spaces $\Hom_{\A}(M,N)$ for all $M,N\in \modd A$]   is called a 
BGG-\textit{algebra}
 if    $\delta(P(i)) \cong I(i)$ for all $i\in Q_0(\A)$ (see \cite[Remark 1.4]{Chan}).
 
  The  functor  $\delta$   for  an   algebra  $A$ in  the previous   examples  is induced by  the anti-automor-phism  defined by  $\epsilon :A\rightarrow A$  via $\epsilon (e_i)=e_i$ and $\epsilon(i\to j)=(j\to i)$ for all $(i\to j)\in Q_1$. 
 The next example shows that 1-quasi-hereditary algebras are in general  not BGG-algebras. 
\\

\textbf{Example 3.} Let  $A_n(C)$  for $n\geq 3$ be the   algebra given by the following   quiver and relations: For all $i,j\in \left\{2, \ldots, n-1\right\}$ the following   holds:

\psset{xunit=1.1mm,yunit=1.1mm,runit=1mm}
\begin{pspicture}(-10,0)(0,0)
\begin{small}
              \rput(0,0){\rnode{0}{$n$}}
              \rput(-13,-10){\rnode{1}{$2$}}
              \rput(-6.5,-10){\rnode{2}{$3$}}
              \rput(-1,-10){\rnode{a}{$\cdots$}}
              \rput(3,-10){\rnode{3}{$i$}}
              \rput(7,-10){\rnode{a}{$\cdots$}}
              \rput(13,-10){\rnode{6}{$n-1$}}
              \rput(0,-20){\rnode{7}{$1$}} 
\end{small}
\begin{footnotesize}
              \rput(73,-10){$
\begin{array}{rcc}
n \ i \ n& = & 0, \\
n\ i\ 1& = & n\ j\ 1 ,\\
1\ i\ n& = & 1\ j\ n  ,\\[5pt]
c_{ij}\cdot \left(i\ n\ j\right)& = &  i\  1\ j,
\end{array}
\ \ \  \  \ \ \ \  C=\left(\begin{array}{ccc}
c_{22}& \cdots & c_{2,n-1}\\
\vdots  & \ddots & \vdots \\
c_{n-1,2}&\cdots & c_{n-1,n-1}
\end{array}\right) \in \GL_{n-2}(K)
$} \end{footnotesize}
              \psset{nodesep=1.5pt,offset=1.5pt,arrows=<-}
              \ncline{0}{1}
              \ncline{1}{0}
              \ncline{0}{2}
              \ncline{2}{0}
              \ncline{0}{3}
              \ncline{3}{0}
              \ncline{0}{6}
              \ncline{6}{0}
              \ncline{1}{7}
              \ncline{7}{1}
              \ncline{2}{7}
              \ncline{7}{2}
              \ncline{3}{7}
              \ncline{7}{3}
              \ncline{6}{7}
              \ncline{7}{6}
\end{pspicture}
\\[2.4cm]
The  order is given by  $1< i< n$  for all $i\in \left\{2,\ldots ,n-1\right\}$, thus  axiom (1) of ~\ref{def1qh} holds.
\\
It's easy to verify that    the set $\ger{B}(i)$  forms a $K$-basis of  $P(i)$ for any  $k\in Q_0$, where
\\[7pt]
\begin{small}
$\begin{array}{l}
\ger{B}(1):=\underbrace{\left\{ e_1, \  (1\ 2\ n \ 2 \ 1) \right\}\cup\{(1\ i\ 1) \mid 2\leq i\leq n-1\}}_{\subset P(1)_1}\cup \{\underbrace{(1\ i),\   (1\  2\  n \   i)}_{\subset P(1)_i} \mid 2\leq i\leq n-1\} \cup \underbrace{\{(1\  2 \  n) \} }_{\subset P(1)_n}\\[6pt]
\ger{B}(i):=
\{\underbrace{(i\ 1),  \  (i\ n \ i\ 1) }_{\subset P(i)_1} ,  \  \underbrace{e_i, \ (i\ n\ i)}_{\subset P(i)_i} \}\cup \{\underbrace{(i\ n \ j) }_{\in  P(i)_j}\mid 2\leq j\leq n-1,\  j\neq i\}\cup \underbrace{\{(i\ n)\}}_{\subset P(i)_n}, \  \      2\leq i\leq n-1, \\[6pt]
\ger{B}(n):= \underbrace{\{( n\ 2\ 1)\}}_{\subset P(n)_1}\cup \{\underbrace{(n\ i)}_{\in  P(n)_i}\mid 2\leq i\leq n-1\}\cup \underbrace{\{e_n\}}_{\subset P(n)_n}.
\end{array}
$
\end{small}        
 \\[5pt]
Axiom (2): The $A$-map $f_{(jk)}:P(k)\to P(j)$ with  $f(e_k)=(j \to k)$  is injective for  $(j,k)\in \left\{(i,n), (1,i)\mid 1< i< n\right\}$.   
 Moreover we have   $\Hom_A(P(n),P(j))=\spann_K\left\{f_{(jn)}\right\}$ and    $\Hom_A(P(i),P(j))=\spann_K\left\{f_{(jn)}\circ f_{(ni)}\right\}$   for all $1< i\neq j < n$.   The definition of $\De(j)$  and   $\im (f_{(jn)}\circ f_{(ni)})\subset \im (f_{(jn)})$  implies    $\De(j)=P(j)/\left(\sum_{j\not\geqslant i}\sum_{f\in \Hom_A(P(i), P(j))}\im (f)\right)=P(j)/(\im (f_{(jn)}))=\spann_K\left\{ \overline{e_j}, \overline{(j 1)}\right\}$ for $1<j<n$ as well as  $\De(1)\cong S(1)$,  $\De(n)=P(n)$.
 For all $i,j\in Q_0 $ with $i\ma j$ we have  $[\De(j):S(i)]=1$.  
    
     The filtration $0\subset M_1\subset M_2\subset  \cdots \subset M_{n-1}\subset P(1)$  with $M_1=\im (f_{(12)}\circ f_{(2n)}) \cong \De(n)$  and  $M_k=\sum_{i=2}^{k}\im  (f_{(1i)})$  is $\De$-good,       since $ P(1)/ M_{n-1}\cong S(1)\cong \De(1)$ and     
$M_k/M_{k-1}\cong \im (f_{(1k)})/\left(\im (f_{(1k)})\cap \sum_{m=2}^{k-1}\im (f_{(1m)})\right)\cong P(k)/\im (f_{(kn)}) \cong \De(k)$  for any  $2\leq k\leq n-1$.
  The filtration $0\subset \im (f_{(jn)})\subset P(j)$   is  $\De$-good  for  any  $1< j < n$. Thus $\left(P(j):\De(i)\right)=1$ for all $i,j\in Q_0$ with $j\ma i$.
\\[5pt]
Axiom (3): Since    $\soc \De(i)\cong S(1)$  for all $i\in Q_0$ and $P(1)\in \ger{F}(\De)$,  we have
 $\soc P(1)\cong S(1)^m$ for  $m\in \NN$. A simple submodule of $P(1)$   is generated by some  non zero  element  $q\in P(1)_1$  with 
$(1 \to j)\cdot q=0$ for all $2\leq j\leq n-1$.    
The basis   $\ger{B}(1)$   of  $P(1)$  shows
$ q=\lambda_1 e_1+ \sum_{i=2}^{n-1}\lambda_i(1\ i\ 1)+ \lambda_n(1\ 2\ n\ 2\ 1)$ with  $\lambda_i\in K$.
  If $\lambda_1\neq 0$, then  $\left\langle q\right)= A\cdot q = P(1)(\neq \soc P(1))$.   Let  $\lambda_1=0$, for every  $j\in \left\{2,\ldots , n-1\right\}$    it is 
\\[5pt]
$
\begin{array}{ccl}
0=(1\to j)\cdot q & = & \displaystyle \sum_{i=2}^{n-1}\lambda_i(1\ i\ 1\ j)+ \lambda_n(1\ 2\ n\   2\ 1\ j) \\
& = & \displaystyle \sum_{i=2}^{n-1}\lambda_i c_{ij}(1\ i\ n\ j)+ \lambda_n c_{2j}(1\ 2\ \underbrace{n\   2\ n}_{=0}\ j)= \displaystyle \sum_{i=2}^{n-1}\lambda_i c_{ij}(1\ j\ n\ j) 
\end{array}
$  
\\
if and only if    \   $ \sum_{i=2}^{n-1}\lambda_i c_{ij}=0$. 
 Since  $\text{det}\:C\neq 0 $, we obtain  $\lambda_i=0$  for all $ 2\leq i\leq n-1$.
 Hence,  $q\in P(1)_1$ generates a simple module  if $q=\lambda (1 2 n 2 1)$ for some  $\lambda\in K\backslash \left\{0\right\}$. 
 For all $\lambda,\mu\in K\backslash \left\{0\right\}$ we have  $\left\langle \lambda (1 2 n 2 1)\right) = \left\langle \mu (1 2 n 2 1)\right)$. Thus 
   $\left\langle 1 2 n 2 1\right) = \soc P(1)\cong S(1)$. Since    $f_{(1i)}:P(i)\hookrightarrow P(1)$,   we obtain  $\soc P(i)\cong S(1)$ for all $i\in Q_0$. The algebra  $\left(A_n(C)\right)^{op}$ is defined by the quiver and relations of   $A_n(C^{\text{tr}})$.  
Using  the same  procedure  we obtain  $\soc P_{A_n(C)^{op}}(i)\cong S(1)$ and  the standard duality implies $\topp I(i)\cong S(1)$ for all $i\in Q_0$.
\\[5pt]
Axiom (4): 	For any $1 < i <n$ we have  $\ker \left(
f_{(ni)}:P(i)\rightarrow  P(n) \right) = \im \left(f_{(in)}:P(n)\hookrightarrow P(i)\right)$, therefore $\De(i)\cong P(i)/\im (f_{(in)})\cong \im (f_{(ni)})\hookrightarrow P(n)=\De(n)$. Since   $\De(1)\cong \soc \De(n)$, we obtain $\De(i)\hookrightarrow \De(n)$, for all $i\in Q_0$. 
Using  the same  procedure  we obtain also  $\De_{A^{op}}(i)\hookrightarrow \De_{A^{op}}(n)$. The standard duality provides $\nabla(n)\twoheadrightarrow \nabla(i)$ for every $i\in Q_0$. The algebra $A_n(C)$ is 1-quasi-hereditary.

For the algebra   $A:= A_4(C)$  with   $C= \begin{scriptsize}\left(
\begin{array}{cc}
1 & q \\
0 & 1
\end{array}
\right) \end{scriptsize}$   and $q \neq 0$ 
    the submodules  of $P(3)$ and  of  $P_{A^{op}}(3)$  are  represented in the    submodule  diagrams (here $\left\langle p\right)=A\cdot p$ resp.   $\left\langle p\right)=A^{op}\cdot p$):
\\ 
\psset{xunit=0.65mm,yunit=0.7mm,runit=1mm}
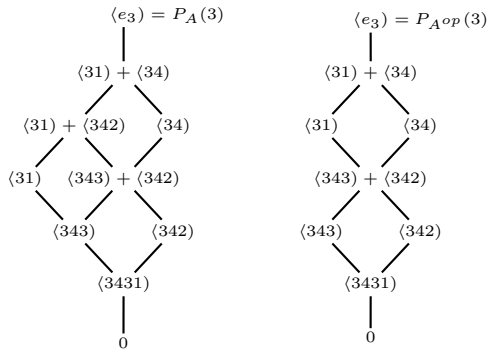
\begin{pspicture}(-23,8)(0,0)
\begin{tiny}
              \rput(0,1){\rnode{0}{$\  \  \ \ \ \ \ \ \ \ \ \  \  \left\langle e_3\right)=P_A(3)$}}
              \rput(0,-10){\rnode{1}{$\left\langle  31\right) + \left\langle  3 4\right)$}}
              \rput(-10,-20){\rnode{2}{$\left\langle 31\right)+\left\langle 342\right)$}}
              \rput(10,-20){\rnode{3}{$\left\langle 34\right)$}}
              \rput(0,-30){\rnode{4}{$\left\langle 343\right)+\left\langle 342\right)$}}
              \rput(-20,-30){\rnode{5}{$\left\langle 31\right)$}}
              \rput(-10,-40){\rnode{6}{$\left\langle 343\right)$}}
              \rput(10,-40){\rnode{7}{$\left\langle 342\right)$}}
              \rput(0,-50){\rnode{8}{$\left\langle 3431\right)$}}
              \rput(0,-61){\rnode{9}{$0$}}
    
   \rput(50,0){ \rput(0,0){\rnode{0x}{$\  \  \  \ \    \ \ \ \ \ \ \ \ \  \  \left\langle e_3\right)=P_{A^{op}}(3)$}}
              \rput(0,-10){\rnode{1x}{$\left\langle  31\right) + \left\langle  34\right)$}}
              \rput(-10,-20){\rnode{2x}{$\left\langle 31\right)$}}
              \rput(10,-20){\rnode{3x}{$\left\langle 34\right)$}}
              \rput(0,-30){\rnode{4x}{$\left\langle 343\right)+\left\langle 342\right)$}}
              \rput(-10,-40){\rnode{6x}{$\left\langle 343\right)$}}
              \rput(10,-40){\rnode{7x}{$\left\langle 342\right)$}}
              \rput(0,-50){\rnode{8x}{$\left\langle 3431\right)$}}
              \rput(0,-60){\rnode{9x}{$0$}}}
              
\end{tiny}
\rput[l](87,-29){\parbox{9cm}{ By duality,  the number of factor modules    of $P_{A^{op}}(3)$ is  the number of submodules   of $I_A(3)$: 
   
$
10=  \left|\left\{\text{submodules of } P(3)\right\}\right| $   
\\
\hspace*{2mm}$
9=  \left|\left\{\text{submodules of } I(3)\right\}\right| $  
\\[5pt]
Therefore on $\modd A$  there can not exist  a duality functor  $\delta$ with $\delta (P(3)) \cong I(3)$. 
\\[3pt]
The  1-quasi-hereditary algebra  $A_4\begin{tiny}\left(
\begin{array}{cc}
1 & q \\
0 & 1
\end{array}
\right)\end{tiny}$  with $q\neq 0$ is not  a  BGG-algebra.}}
              \psset{nodesep=1.5pt}
              \ncline{0}{1}
              \ncline{1}{2}
              \ncline{1}{3}
              \ncline{2}{4}
              \ncline{2}{5}
              \ncline{3}{4}
              \ncline{4}{6}
              \ncline{5}{6}
              \ncline{4}{7}
              \ncline{6}{8}
              \ncline{7}{8}
              \ncline{8}{9}
              \ncline{0x}{1x}
              \ncline{1x}{2x}
              \ncline{1x}{3x}
              \ncline{2x}{4x}
              \ncline{3x}{4x}
              \ncline{4x}{6x}
              \ncline{4x}{7x}
              \ncline{6x}{8x}
              \ncline{7x}{8x}
              \ncline{8x}{9x}
\end{pspicture}
\\[4.6cm]
If  $q=0$ then $A_4(C)$  is  a  BGG-algebra  with the duality   induced by  an anti-automorphism.
\\

\begin{small}Specific for the class of quasi-hereditary algebras is the concept of Ringel-duality: 
Let $\A$ be a quasi-hereditary algebra. For any $i\in Q_0(\A)$ there exists up to isomorphism an unique indecomposable  $\A$-module $T(i)$ having a $\De$-good and $\nabla$-good filtration with $ (T(i): \De(i))=(T(i):\nabla(i))= [T(i):S(i)]= 1$  and   $(T(i): \De(j))=(T(i):\nabla(j))=[T(i):S(j)]= 0$ for all $j\not\ma i$,  moreover there exists   a submodule $Y(i)\in \ger{F}(\nabla)$  of $T(i)$ with $T(i)/Y(i)\cong \nabla(i)$   
and  a factor module $X(i)\in \ger{F}(\De)$  with  $\Kern (T(i)\twoheadrightarrow X(i)) \cong \De(i)$. 
The algebra $R(\A):=\End_{\A}(\bigoplus_{i\in Q_0}T(i))^{op}$ with the opposite order $\geqslant$    is also quasi-hereditary, where $T=\bigoplus_{i\in Q_0}T(i)$ is the \textit{characteristic tilting module} of $\A$.   In particular,  $\ger{F}(\De)\cap \ger{F}(\nabla) = \add (T)$  and $\A $  is isomorphic to  $ R(R(\A))$  as  a  quasi-hereditary algebra  (see  \cite{Rin1}).  The algebra $R(\A)$ is called \textit{Ringel-dual} of $\A$. \end{small}
\\

Some properties of the characteristic tilting module and the Ringel-dual of a 1-quasi-hereditary algebra  are considered  in section 5 and 6 in \cite{P}.  
According to   \cite[Remark 5.3]{P} for the direct summands of the  characteristic tilting $A_n(C)$-module      we obtain   $T(1)\cong S(1)$, $T(n)\cong P(1)$ and  $T(i) \cong P(1)/\left(\sum_{j=2\atop j\neq i}^{n-1}P(j)\right)\cong \bigcap_{j=2\atop j\neq i}^{n-1}\ker\left(P(1)\twoheadrightarrow I(j)\right)$ for all $2\leq i\leq n-1$,     since any vertex $i\in Q_0\backslash\left\{1,n\right\}$ is a 	neighbor  of $1$.   Consequently,  $T(i)$ is  a submodule (and a factor module) of $P(1)$, thus  there exists some   element in $ P(1)_1$, which generates $T(i)$. In particular,  the Ringel-dual of $A_n(C)$ is also 1-quasi-hereditary  (see  \cite[Theorem 6.1]{P}).

\begin{num}\begin{normalfont}\textbf{Lemma.}\end{normalfont} Let $C=(c_{ij})_{2\leq i,j\leq n-1}\in \GL_{n-2}(K)$,  the algebra $A_n(C)$ given in Example 3 and $C^{-1}=(d_{ij})_{2\leq i,j\leq n-1}$. Then  the following hold:
\begin{itemize}
	\item[(1)] $\ger{t}(i):= \displaystyle \sum_{j=2}^{n-1}d_{ij}\cdot (1\ j \ 1)$ generates $T(j)$ for $2\leq i\leq n-1$.  \\ 
	Moreover, $\ger{t}(1):=(1  2n  2 1)$  generates $T(1)$ and   $\ger{t}(n):=e_1$ generates  $T(n)$. 
	\item[(2)] $R(A_n(C))\cong A_n(C^{-1})$.
\end{itemize}
\end{num}

\textit{Proof.} \textit{(1)}  The $A_n(C)$-module $\left\langle \ger{t}(i)\right)$ is local with $\topp \left\langle \ger{t}(i)\right) \cong S(1)\cong \De(1)$ since  $\ger{t}(i)\in P(1)_1$ for all $i\in Q_0$. Using the calculations in  Example 3,  we obtain $T(1)\cong \left\langle 12n21\right)$ and $T(n)\cong \left\langle e_1\right)$, since $T(1)\cong \soc P(1)$ and $T(n)\cong P(1)$. For all $i,k\in \left\{2,\ldots ,n-1\right\}$ we have 
\begin{center}
$
\begin{array}{ccl}
(1\to k)\cdot \ger{t}(i)& = & \displaystyle \sum_{j=2}^{n-1}d_{ij}\cdot (1\ j\ 1 \ k) = \sum_{j=2}^{n-1}d_{ij}\cdot c_{jk} \cdot (1\ j \ n \ k) \\
&= & \displaystyle \left(\sum_{j=2}^{n-1}d_{ij}\cdot c_{jk} \right)\cdot (1 \ k \ n\ k) = \begin{small}\left\{
\begin{array}{cl}
(1 \ i \ n \ i) & \text{if } k=i,\\
0 & \text{else.}
\end{array}\right.
 \end{small}
\end{array}$
\end{center}
 Consequently,  
    $\rad \left\langle \ger{t}(i)\right) =\left\langle  1 \ i \ n \ i \right) \cong \De(i)$ and $0\subset \soc \left\langle \ger{t}(i)\right) = \left\langle \ger{t}(1)\right) \subset \rad \left\langle \ger{t}(i)\right) \subset  \left\langle \ger{t}(i)\right)  $ is the unique Jordan-Hölder-filtration of $\left\langle \ger{t}(i)\right)$.  The filtration $0\subset \rad \left\langle \ger{t}(i)\right) \subset  \left\langle \ger{t}(i)\right) $ resp. $0\subset \soc \left\langle \ger{t}(i)\right) \subset \left\langle \ger{t}(i)\right)$ is $\De$-good resp. $\nabla$-good with the properties of $T(i)$, thus $\left\langle \ger{t}(i)\right)\cong T(i)$  for every  $2\leq i\leq n-1$.

\textit{(2)} We consider $T(i)$ as a submodule of $P(1)$.  Since $R(A_n(C))\cong \End_{A_n(C)}(T)^{op}$  is 1-quasi-hereditary, the quivers  of $(A_n(C),\ma)$ and $(R(A_n(C)),\geqslant)$ have  the same shape (see  \cite[Theorem 2.7]{P}).  
The vertex $i$ in   $Q_0(R(A_n(C)))$   corresponds to the direct summand $T(i)$ of $T= \bigoplus_{i\in Q_0}T(i)$. 
 For $(l,m)\in \left\{(1,j), (j,1),  (n,j), (j,n) \mid 2\leq j\leq n-1\right\}$ 
 we denote  by $\tau_{(l,m)}$  the following maps: 
\\[3pt]
\hspace*{-2mm}$
\begin{array}{lll}
\tau_{(1,j)}:T(1)\to T(j) \text{ with   }  \tau_{(1,j)}(\ger{t}(1))=\ger{t}(1) &\text{and}& \tau_{(j,1)}:T(j)\to T(1) \text{ with  }  \tau_{(j,1)}(\ger{t}(j))=\ger{t}(1)\\
\tau_{(j,n)}:T(j)\to T(n)  \text{ with   } \tau_{(j,n)}(\ger{t}(j))=\ger{t}(j)&\text{and}& \tau_{(n,j)}:T(n)\to T(j) \text{ with   } \tau_{(n,j)}(\ger{t}(n))= \ger{t}(j).
\end{array}
$
\\[3pt] 
It is easy to compute that 
 the space of maps  in  $\Hom_{A_n(C)}(T(l),T(m))$ which  do not   factors through  $\add \left(T\right)$  is spanned  by $\tau_{(l,m)}$, thus  $\tau_{(l,m)}$ corresponds to the arrow $(l\to m)$ for  any $(l,m)$.  We can compute the relations:
For all $2\leq i,j\leq n-1$ we have $\tau_{(i,n)}\circ \tau_{(1,i)} = \tau_{(j,n)}\circ\tau_{(1,j)}$ and $\tau_{(i,1)}\circ \tau_{(n,i)} = \tau_{(j,1)}\circ\tau_{(n,j)}$, thus $(1\ j\ n)=(1\ i\ n)$ and $(n \ j\ 1)=(n\ i \  1)$. 
 Moreover,    $\left(\tau_{(1,j)}\circ \tau_{(i,1)}\right)(\ger{t}(i))=\ger{t}(1)$ and $\left(\tau_{(n,j)}\circ \tau_{(i,n)}\right)(\ger{t}(i))=\ger{t}(i)\cdot \ger{t}(j) $ with 
\begin{center}
$\begin{array}{ccl}
\ger{t}(i)\cdot \ger{t}(j)	& = &\displaystyle  \left(\sum_{l=2}^{n-1}d_{il}\cdot (1\ l\ 1)\right) \cdot \left(\sum_{k=2}^{n-1}d_{jk}\cdot (1\ k\ 1)\right) = \displaystyle \sum_{l=2}^{n-1}d_{il}\sum_{k=2}^{n-1} d_{jk}\cdot (1\ k\ 1\ l\ 1)\\
& = & \displaystyle \sum_{l=2}^{n-1}d_{il}\sum_{k=2}^{n-1}\left( d_{jk}\cdot c_{kl}\right) \cdot (1 \ k \ n \ l  \ 1)  =   d_{ij}\cdot  \ger{t}(1)
\end{array}
$
\end{center}
We obtain  $\tau_{(n,j)}\circ\tau_{(i,n)} =d_{ij}\left( \tau_{(1,j)} \circ \tau_{(i,1)}\right)$, 
 thus   $(i\ n\ j)= d_{ij} \cdot (i \ 1\ j) $. The map $\tau_{(j,1)}\circ\tau_{(1,j)}$ 
  is  zero-map, thus $(1\ j\ 1)=0$ for any $j\in Q_0\backslash\left\{1,n\right\}$. There are  no relations between  paths  starting and ending in $n$, since  the set  $\left\{\left(\tau_{(j,n)}\circ \tau_{(n,j)}\right)(\ger{t}(n))=\ger{t}(j) \mid 2\leq j\leq n-1\right\} \cup \left\{\ger{t}(1)\right\}$ is linearly independent.

By  interchanging the notations $1\mapsto n $ and $n\mapsto 1$ we obtain  that  the quiver and relations of $R(A_n(C))$ are those  of the algebra $A_n(C^{-1})$. \hfill $\Box$
\\ 

The class of 1-quasi-hereditary algebras is not closed under Ringel-duality. According to   \cite[Theorem 6.1]{P} the algebra  $R(A)$ is 1-quasi-hereditary if and only if the factor algebra $A(i)= A/\left(A\left(\sum_{j\not\ma i}e_j\right)A\right)$ is 1-quasi-hereditary for every $i\in Q_0$ (see  \cite[Section 5]{P}). The next example presents a 1-quasi-hereditary algebra $A$ such that  there exists $i\in Q_0(A)$ with 
  $A(i)$  not being   1-quasi-hereditary.
\\

\textbf{Example 4.} The algebra $A$ given by the following quiver and relations      is 1-quasi-hereditary with the partial order   $1< 2< 3< 5< 6$ and  $1< 4< 5$. 
 \\[-8pt]

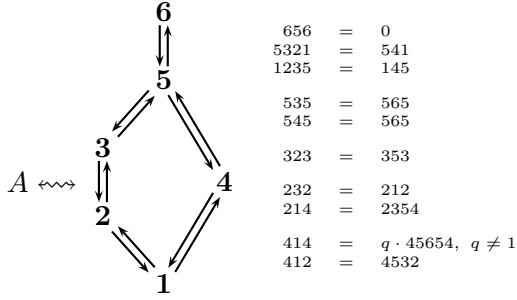
\psset{xunit=0.8mm,yunit=0.9mm,runit=0.8mm}
\begin{pspicture}(-20,0)(110,0)
\begin{small}
\rput(-20,-25){\rnode{A}{$A \leftrightsquigarrow $}}
              \rput(0,0){\rnode{0v}{\textbf{6}}}
              \rput(0,-10){\rnode{1v}{\textbf{5}}}
              \rput(-10,-20){\rnode{2v}{\textbf{3}}}
              \rput(10,-25){\rnode{3v}{\textbf{4}}}
              \rput(-10,-30){\rnode{4v}{\textbf{2}}}
              \rput(0,-40){\rnode{5v}{\textbf{1}}} 
              
\end{small}
\begin{tiny}
\rput(38,-20){$
\begin{array}{rcl}
656&=& 0 \\
5321&=&541 \\
1235&=&145\\[6pt]
535&=&565 \\
545&=& 565\\[6pt]
323&=&353 \\[6pt]
232&=&212 \\
214&=&2354 \\[6pt]
414&=&  q\cdot  45654 ,  \   q\neq 1\\
412&=&4532 \\
\end{array}
$}
\end{tiny}
\rput(120,-20){\parbox{9cm}{\begin{footnotesize}
The Auslander-algebra $\A_m$ of $K[x]/\left\langle x^m\right\rangle$  is  obtained   by  the construction in Example 2 with $B=K[x]/\left\langle x^m\right\rangle$ and the set of ideals $\ger{X}=\left\{\mathcal{X}(i):=\left\langle \overline{x}^{i}\right\rangle\mid  1\leq i\leq m \right\}$.
\\
For  $i\in Q_0 $  with     $ i=1,2,3$  the algebra  $A(i)$ is isomorph to the Auslander-algebra  $\A_{i}$     and  $A(4)\cong \A_2$. 

Thus $A(i)$  is 1-quasi-hereditary for   $i=1,2,3,4,6$   and    we can compute    an element $\ger{t}(i)\in P(1)_1$   which   generates     $ T(i)$:
\\[5pt]
$
\begin{array}{lll}
 \ger{t}(1) = (1456541), \hspace{0.1cm} & \ger{t}(3)=  q \cdot (12321)-(141), & \hspace{0.1cm} \ger{t}(6) =e_1.\\
\ger{t}(2)=  (14541),& \ger{t}(4)=   (141)-(12321). &  
\end{array}
$
 \end{footnotesize}}}
\psset{nodesep=1.5pt,offset=1.5pt,arrows=<-}
              \ncline{0v}{1v}
              \ncline{1v}{0v}
              \ncline{1v}{2v}
              \ncline{2v}{1v}
              \ncline{1v}{3v}
              \ncline{3v}{1v}
              \ncline{2v}{4v}
              \ncline{4v}{2v}
              \ncline{4v}{5v}
              \ncline{5v}{4v}
              \ncline{3v}{5v}
              \ncline{5v}{3v}
\end{pspicture}
\\[3cm]

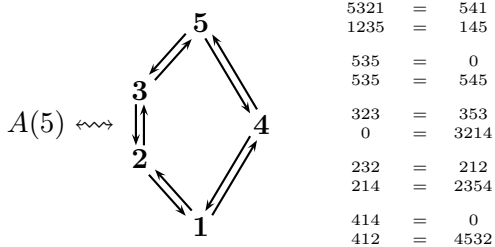
\psset{xunit=0.8mm,yunit=0.9mm,runit=0.8mm}
\begin{pspicture}(-20,0)(110,0)
\begin{small}
\rput(-23,-25){\rnode{A}{$A(5) \leftrightsquigarrow $}}              \rput(0,-10){\rnode{1}{\color{black}{\textbf{5}}}}
              \rput(-10,-20){\rnode{2}{\color{black}{\textbf{3}}}}
              \rput(10,-25){\rnode{3}{\color{black}{\textbf{4}}}}
              \rput(-10,-30){\rnode{4}{\color{black}{\textbf{2}}}}
              \rput(0,-40){\rnode{5}{\color{black}{\textbf{1}}}} 
              
\end{small}
\begin{tiny}
\rput(36,-25){$
\begin{array}{ccc}
5321&=&541 \\
1235&=&145\\[6pt]
535&=&0 \\
535&=& 545\\[6pt]
323&=&353 \\
0&=&3214 \\[6pt]
232&=&212 \\
214&=&2354 \\[6pt]
414&=&  0\\
412&=&4532 \\
\end{array}
$}
\end{tiny}
\rput(120,-26){\parbox{9.2cm}{\begin{footnotesize}
The elements 
$p_1:=(14541)$ and $p_2:=(141)-(12321)$ in $P_{A(5)}(1)_1$ are linear independent with   $(12)\cdot p_1 =(14)\cdot p_1=0$ and $(12)\cdot p_2 =(14)\cdot p_2=0$. The submodules of $P_{A(5)}(1)$ generated by $p_1$ and $p_2$ are simples. Hence $\soc P_{A(5)}(1) \not\cong S(1)$,  thus   the algebra  $A(5)$ is not 1-quasi-hereditary.
\\[5pt]
The $A$-module   $T(5)$   is a submodule  of  $P(1)\oplus P(1)$   generated by $\big((121),((141)-(12321))\big)$ and $\big((141),0\big)$.
 \end{footnotesize}}}
\psset{nodesep=1.5pt,offset=1.5pt,arrows=<-}
              \ncline{0}{1}
              \ncline{1}{0}
              \ncline{1}{2}
              \ncline{2}{1}
              \ncline{1}{3}
              \ncline{3}{1}
              \ncline{2}{4}
              \ncline{4}{2}
              \ncline{4}{5}
              \ncline{5}{4}
              \ncline{3}{5}
              \ncline{5}{3}
\end{pspicture}
\\[4cm]

The quiver of a  1-quasi-hereditary  algebra  $\A$ depends only on  $(Q_0(\A),\ma)$: If  $i\stackrel{\alpha}{\rightarrow}j\in Q_1(\A)$, then 
 $i\triangleleft j$ resp.  $i\triangleright j$  and   $\left|\left\{\alpha\in Q_1(\A)\mid i\stackrel{\alpha}{\rightarrow}j\right\}\right|= \left|\left\{\alpha\in Q_1(\A)\mid j\stackrel{\alpha}{\rightarrow}i\right\}\right|=1$  (see  \cite[Theorem  2.7]{P}). All Jordan-Hölder-filtrations of $\De(i)$ and $\nabla(i)$ as well as good filtrations of $P(i)$ and $I(i)$ are connected
 to  the sequences  from $\mathcal{T}(i)$ resp. $\mathcal{L}(i)$ which 
also depend only  on the structure of $\ma$ for every $i\in Q_0$ (see \cite[Propositions 4.1 and 4.2]{P}). Each   of these  sequences can be completed to some sequence  from $\mathcal{T}(n)=\mathcal{L}(1)=\left\{\textbf{i}=(i_1, i_2, \ldots , i_n)\in Q_0(\A)^{n} \mid   i_{k}\not\geqslant i_{t}, \ 1\leq  k < t \leq n \right\}$, thus we have 
\\[5pt]
$\mathcal{L}(1)\leftrightarrow \begin{footnotesize}\left\{
\begin{array}{c}
\text{Jordan-Hölder}\\
\text{filtr. of  }\De(n)	\\[4pt]
 \mathcal{S}(\textbf{i})
\end{array}
\right\} \end{footnotesize}\leftrightarrow \begin{footnotesize}\left\{
\begin{array}{c}
\text{Jordan-Hölder.}\\
\text{filtr. of  }\nabla(n)	\\[4pt]
 \widetilde{\mathcal{S}}(\textbf{i})
\end{array}
\right\} \end{footnotesize} \leftrightarrow \begin{footnotesize}\left\{
\begin{array}{c}
\De\text{-good filtr.}\\
\text{of  }P(1)	\\[4pt]
 \mathcal{D}(\textbf{i})
\end{array}
\right\} \end{footnotesize} \leftrightarrow \begin{footnotesize}\left\{
\begin{array}{c}
\nabla\text{-good filtr.}\\
\text{of  }I(1)	
\\[4pt]
 \mathcal{N}(\textbf{i})\end{array}
\right\} \end{footnotesize}$
\\[5pt]
The sequence of indices of the simple factors  in $\widetilde{\mathcal{S}}(\textbf{i})$ and $\De$-good  factors in $\mathcal{D}(\textbf{i})$ are the same (and in the same order). The same holds for the indices of factors in  $\mathcal{S}(\textbf{i})$ and $\mathcal{N}(\textbf{i})$. 	Furthermore the indices of simple factors of $\widetilde{\mathcal{S}}(\textbf{i})$ and $\mathcal{S}(\textbf{i})$ are the same, but in reversed  order.  Since $\De(n)$ (resp. $\nabla(n)$) has finitely many submodules, all filtrations of $\De(n)$ can by represented in the submodule diagram of $\De(n)$, i.e. the Hasse diagram of $\left(\left\{\text{submodules of } \De(n)\right\}, \subseteq \right)$.

  The last  observation   shows  that   $\nabla$-good filtrations of $I(1)$  can be represented in a diagram (we  will   call   it   $\nabla$-good diagram of $I(1)$)  whose   shape  coincides    with the submodule diagram of $\De(n)$.   By standard  duality,   the submodule diagram of $\nabla(n)$ and the  $\De$-good diagram of $P(1)$ have  the same shape.  Moreover, the shape of 
the diagram  of $\De(n)$ is obtained from the diagram  of $\nabla(n)$ by turning it upside down. In particular,  all  $\nabla$-good filtrations of $I(i)$  
are parts of $\nabla$-diagram  of $I(1)$  resp.  
 the  $\De$-good diagram of  $P(i)$  is a $\De$-subdiagram of $P(1)$ for any $i\in Q_0$.
We illustrate this using  Example 4.
	 
 The set $\mathcal{L}(1)=\left\{\left(1,2,3,4,5,6\right), \  \left(1,2,4,3,5,6\right), \  \left(1,4,2,3,5,6\right)\right\}$ shows that there exist three  Jordan-Hölder filtrations of $\De(6)$ resp. $\nabla(6)$  and three    good filtrations of $P(1)$ resp.  $I(1)$. All  Jordan-Hölder filtrations of $\De(j)$  resp.  $\nabla(j)$  are subfiltrations of $\De(6)$ resp. factor filtrations of $\nabla(6)$ (in the picture $\ger{K}(j)=\ker (\nabla(n)\twoheadrightarrow \nabla(j))$. Similar   $\De$-good filtrations of $P(j)$  are  subfiltrations of  $P(1)$ and  $\nabla$-good filtrations of $I(1)$ are  factor filtrations of those  of  $I(1)$ (here $\mathcal{K}(j)=\ker(I(1)\twoheadrightarrow I(j))$) for all $j\in Q_0$.   
\\
\\

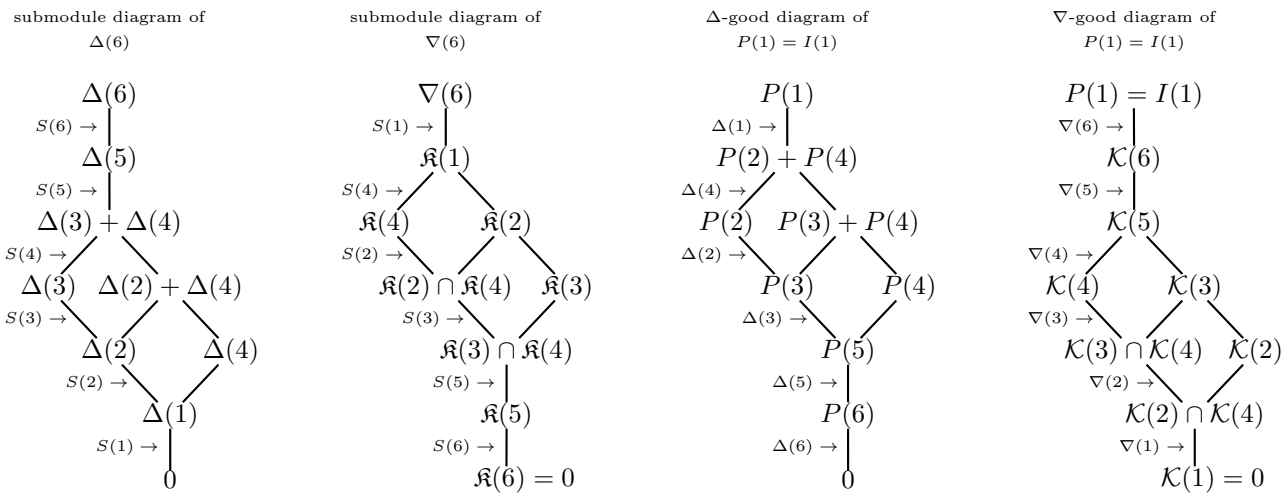
\psset{xunit=0.8mm,yunit=0.85mm,runit=1mm}
\begin{pspicture}(50,0)(20,0)
\rput(25,0){\begin{footnotesize}
\rput(0,0){\rnode{1}{$\De(6)$}}
\rput(0,-10){\rnode{2}{$\De(5)$}}
\rput(0,-20){\rnode{3}{$\De(3)+\De(4)$}}
\rput(-10,-30){\rnode{4}{$\De(3)$}}
\rput(10,-30){\rnode{5}{$\De(2)+\De(4)$}}
\rput(0,-40){\rnode{6}{$\De(2)$}}
\rput(20,-40){\rnode{7}{$\De(4)$}}
\rput(10,-50){\rnode{8}{$\De(1)$}}
\rput(10,-60){\rnode{9}{$0$}}
\end{footnotesize}
\begin{tiny}
\rput(0,12){\rnode{s}{submodule diagram of}}
\rput(0,8){\rnode{ss}{$\De(6)$}}
\rput(-7,-5){\rnode{1x}{$S (6) \rightarrow $}}
\rput(-7,-15){\rnode{23x}{$S (5) \rightarrow$}}
\rput(-12,-25){\rnode{3x}{$S (4) \rightarrow$}}
\rput(-12,-35){\rnode{45x}{$S (3) \rightarrow$}}
\rput(-2,-45){\rnode{5x}{$S (2)\rightarrow $}}
\rput(3,-55){\rnode{6x}{$S (1)\rightarrow$}}
             \end{tiny}}
             \ncline{1}{2}
             \ncline{2}{3}
             \ncline{3}{4}
             \ncline{3}{5}
             \ncline{4}{6}
             \ncline{5}{6}
             \ncline{5}{7}
             \ncline{6}{8}
             \ncline{7}{8}
             \ncline{8}{9}
             
\end{pspicture} 
\psset{xunit=0.8mm,yunit=0.85mm,runit=1mm}
\begin{pspicture}(50,0)(-5,0)
\rput(35,-60){\begin{footnotesize}
\rput(0,0){\rnode{1}{$ \ \ \ \ \ger{K}(6)=0$}}
\rput(0,10){\rnode{2}{$\ger{K}(5)$}}
\rput(0,20){\rnode{3}{$\ger{K}(3)\cap \ger{K}(4)$}}
\rput(10,30){\rnode{4}{$\ger{K}(3)$}}
\rput(-10,30){\rnode{5}{$\ger{K}(2)\cap \ger{K}(4)$}}
\rput(0,40){\rnode{6}{$\ger{K}(2)$}}
\rput(-20,40){\rnode{7}{$\ger{K}(4)$}}
\rput(-10,50){\rnode{8}{$\ger{K}(1)$}}
\rput(-10,60){\rnode{9}{$\nabla(6)$}}
\end{footnotesize}  
\begin{tiny}
\rput(-10,72){\rnode{s}{submodule diagram of}}
\rput(-10,68){\rnode{ss}{$\nabla(6)$}}
\rput(-17,55){\rnode{1x}{$S (1) \rightarrow $}}
\rput(-22,45){\rnode{23x}{$S (4) \rightarrow$}}
\rput(-22,35){\rnode{3x}{$S (2) \rightarrow$}}
\rput(-12,25){\rnode{45x}{$S (3) \rightarrow$}}
\rput(-7,15){\rnode{5x}{$S (5)\rightarrow $}}
\rput(-7,5){\rnode{6x}{$S (6)\rightarrow$}}
             \end{tiny} }
             \ncline{1}{2}
             \ncline{2}{3}
             \ncline{3}{4}
             \ncline{3}{5}
             \ncline{4}{6}
             \ncline{5}{6}
             \ncline{5}{7}
             \ncline{6}{8}
             \ncline{7}{8}
             \ncline{8}{9} 
             
\end{pspicture} 
\psset{xunit=0.8mm,yunit=0.85mm,runit=1mm}
\begin{pspicture}(50,0)(-5,0)
\rput(35,-60){\begin{footnotesize}
\rput(0,0){\rnode{1}{$0$}}
\rput(0,10){\rnode{2}{$P(6)$}}
\rput(0,20){\rnode{3}{$P(5)$}}
\rput(10,30){\rnode{4}{$P(4)$}}
\rput(-10,30){\rnode{5}{$P(3)$}}
\rput(0,40){\rnode{6}{$P(3)+P(4)$}}
\rput(-20,40){\rnode{7}{$P(2)$}}
\rput(-10,50){\rnode{8}{$P(2)+P(4)$}}
\rput(-10,60){\rnode{9}{$P(1)$}}
\end{footnotesize}  
\begin{tiny}
\rput(-10,72){\rnode{s}{$\De$-good  diagram of}}
\rput(-10,68){\rnode{ss}{$P(1)=I(1)$}}
\rput(-17,55){\rnode{1x}{$\De (1) \rightarrow $}}
\rput(-22,45){\rnode{23x}{$\De (4) \rightarrow$}}
\rput(-22,35){\rnode{3x}{$\De (2) \rightarrow$}}
\rput(-12,25){\rnode{45x}{$\De (3) \rightarrow$}}
\rput(-7,15){\rnode{5x}{$\De (5)\rightarrow $}}
\rput(-7,5){\rnode{6x}{$\De (6)\rightarrow$}}
             \end{tiny} }
             \ncline{1}{2}
             \ncline{2}{3}
             \ncline{3}{4}
             \ncline{3}{5}
             \ncline{4}{6}
             \ncline{5}{6}
             \ncline{5}{7}
             \ncline{6}{8}
             \ncline{7}{8}
             \ncline{8}{9} 
             
\end{pspicture} 
\psset{xunit=0.8mm,yunit=0.85mm,runit=1mm}
\begin{pspicture}(50,0)(-5,0)
\rput(25,0){\begin{footnotesize}
\rput(0,0){\rnode{1}{$P(1)=I(1)$}}
\rput(0,-10){\rnode{2}{$\mathcal{K}(6)$}}
\rput(0,-20){\rnode{3}{$\mathcal{K}(5)$}}
\rput(-10,-30){\rnode{4}{$\mathcal{K}(4)$}}
\rput(10,-30){\rnode{5}{$\mathcal{K}(3)$}}
\rput(0,-40){\rnode{6}{$\mathcal{K}(3)\cap \mathcal{K}(4)$}}
\rput(20,-40){\rnode{7}{$\mathcal{K}(2)$}}
\rput(10,-50){\rnode{8}{$\mathcal{K}(2)\cap \mathcal{K}(4)$}}
\rput(10,-60){\rnode{9}{$\ \ \ \  \mathcal{K}(1)=0$}}
\end{footnotesize}
\begin{tiny}
\rput(0,12){\rnode{s}{$\nabla$-good diagram of}}
\rput(0,8){\rnode{ss}{$P(1)=I(1)$}}
\rput(-7,-5){\rnode{1x}{$\nabla (6) \rightarrow $}}
\rput(-7,-15){\rnode{23x}{$\nabla (5) \rightarrow$}}
\rput(-12,-25){\rnode{3x}{$\nabla (4) \rightarrow$}}
\rput(-12,-35){\rnode{45x}{$\nabla (3) \rightarrow$}}
\rput(-2,-45){\rnode{5x}{$\nabla (2)\rightarrow $}}
\rput(3,-55){\rnode{6x}{$\nabla (1)\rightarrow$}}
             \end{tiny}}
             \ncline{1}{2}
             \ncline{2}{3}
             \ncline{3}{4}
             \ncline{3}{5}
             \ncline{4}{6}
             \ncline{5}{6}
             \ncline{5}{7}
             \ncline{6}{8}
             \ncline{7}{8}
             \ncline{8}{9}
             
\end{pspicture} \\[5.4cm]
The subquotients are illustrated as follows: The module $N$   pointing   to the line which connects the modules $M$ and $M'$ such that $M\subset M'$ has the meaning $M'/M\cong N$.

\section{Modules generated by the paths  $p(j,i,k)$}
\begin{small} In this section $(A,\ma)$ is a 1-quasi-hereditary algebra. 
In the following we use some notations introduced in  \cite{P}: For any $i\in Q_0$ we have the sets 
\begin{center}
$\Lambda_{(i)}:= \left\{j\in Q_0\mid j\ma i\right\}$\   and \   $\Lambda^{(i)}:=\left\{j\in Q_0\mid j\geqslant i\right\}$.
\end{center}
 If $i$ is a small resp. large  neighbor of $j$ w.r.t. $\ma$,  we write $i\triangleleft j $ resp. $i\triangleright j$.  Recall that $\left\langle p\right)=A\cdot p$.
\end{small}
\\  

We recall the definition of $p(j,i,k)$ from  \cite[Sec. 3]{P}:  Let $j,i,k $ be in  $ Q_0$ with $i\in \Lambda^{(j)}\cap \Lambda^{(k)}$. 
 There  exist    $\lambda_0,\ldots , \lambda_r \in \Lambda^{(j)}$    and $\mu_0, \ldots , \mu_m \in \Lambda^{(k)}$  with   $j=\lambda_0 \ma \lambda_1 \ma \cdots \ma \lambda _r=i$   and  $ i =\mu_0\geqslant \mu_1\geqslant \cdots \geqslant \mu_m=k$  giving  a   path 
\\
  \parbox{13.5cm}{\begin{center}
 $p(j,i,k):= \left(j  \rightarrow \lambda_1 \rightarrow  \cdots  \rightarrow \lambda_{r-1}\rightarrow i\rightarrow \mu_1 \rightarrow \cdots \rightarrow \mu_{m-1}\rightarrow k\right)$
\end{center}
(a path $p(j,i,k)$ [black] resp. $p(k,i,j)$ [gray] is visualised   in the picture to the right).  We fix a path of the form   $p(j,i,i)$  and $p(i,i,k)$,  we obviously have  $p(j,i,k)= p(i,i,k)\cdot p(j,i,i)$. We   denoted by $f_{(j,i,k)}$ the $A$-map corresponding to $p(j,i,k)$, i.e. $f_{(j,i,k)}:P(k)\to P(j)$ with $f_{(j,i,k)}(a\cdot e_k)= a\cdot p(j,i,k)$ for all $a\in A$.
   We obtain   $f_{(j,i,k)}=f_{(j,i,i)}\circ f_{(i,i,k)}$.
   A   path $p(j,i,i)$  resp. $p(i,i,k)$  is 
 }
 \psset{xunit=0.64mm,yunit=0.69mm,runit=1mm}
\begin{pspicture}(0,0)(0,0)
\rput(25,22){
\begin{tiny}
\rput(0,10){\rnode{A}{$n$}}
\rput(0,-44){\rnode{B}{$1$}}
\rput(20,-17){\rnode{aa}{}}
\rput(-20,-17){\rnode{bb}{}}
\rput(-11,3){\rnode{A1}{}}
\rput(-11,-37){\rnode{B1}{}}
\rput(11,3){\rnode{A2}{}}
\rput(11,-37){\rnode{B2}{}}
              \rput(0,1){\rnode{00}{$   \lambda_r =i=\mu_0$}}
              \rput(-11,-32){\rnode{1}{$ \ \ \ \ \  \ \ \  j=\lambda_0$}}
              \rput(-5,-23){\rnode{11}{$\lambda_1$}}
              \rput(-10,-14){\rnode{12}{$\lambda_2$}}
              \rput(-6,-8){\rnode{13}{$\lambda_{r-1}$}}
              \rput(8,-8){\rnode{21}{$\mu_2$}}
              \rput(7,-16){\rnode{22}{$\mu_{m-1}$}}
              \rput(12,-26){\rnode{2}{$\ \ \   \ \ \ k=\mu_m$}}
              \end{tiny}
                      }
             \psset{nodesep=1pt,offset=1.3pt}
              \ncline{->}{1}{11}
              \ncline[linewidth=0.7pt,linecolor=black!70]{->}{11}{1}
              \ncline{->}{11}{12}
              \ncline[linewidth=0.7pt,linecolor=black!70]{->}{12}{11}
              \ncline[linestyle=dotted]{12}{13}
              \ncline{->}{13}{00}
              \ncline[linewidth=0.7pt,linecolor=black!70]{->}{00}{13}
              \ncline{->}{00}{21}
              \ncline[linewidth=0.7pt,linecolor=black!70]{->}{21}{00}
              \ncline{->}{22}{2}
              \ncline[linewidth=0.7pt,linecolor=black!70]{->}{2}{22}
              \ncline[linestyle=dotted]{21}{22}
              \ncline[linewidth=0.7pt,linecolor=black!70]{->}{A}{A1}
              \ncline[linewidth=0.7pt,linecolor=black!70]{->}{A1}{A}
              \ncline[linewidth=0.7pt,linecolor=black!70]{->}{B1}{B}
              \ncline[linewidth=0.7pt,linecolor=black!70]{->}{B}{B1}
              \ncline[linewidth=0.7pt,linecolor=black!70]{->}{A}{A2}
              \ncline[linewidth=0.7pt,linecolor=black!70]{->}{A2}{A}
              \ncline[linewidth=0.7pt,linecolor=black!70]{->}{B2}{B}
              \ncline[linewidth=0.7pt,linecolor=black!70]{->}{B}{B2}
              \ncarc[arcangle=25,linestyle=dotted,linecolor=black!70]{A2}{aa}
              \ncarc[arcangle=-25,linestyle=dotted,linecolor=black!70]{->}{aa}{A2}
              \ncarc[arcangle=25,linestyle=dotted,linecolor=black!70]{->}{aa}{B2}
              \ncarc[arcangle=-25,linestyle=dotted,linecolor=black!70]{B2}{aa}
              \ncarc[arcangle=-25,linestyle=dotted,linecolor=black!70]{A1}{bb}
              \ncarc[arcangle=-25,linestyle=dotted,linecolor=black!70]{->}{bb}{B1}
              \ncarc[arcangle=25,linestyle=dotted,linecolor=black!70]{B1}{bb}
              \ncarc[arcangle=25,linestyle=dotted,linecolor=black!70]{->}{bb}{A1}
\end{pspicture} 
\\
	increasing  
resp.  declining   and $e_i=p(i,i,i)$  is the  trivial path $e_j$. In the following,  we  point out of some  facts which are presented in  \cite{P}:
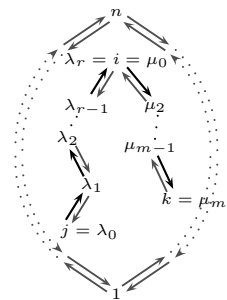
\begin{itemize}
\item[(1)] The set $\left\{p(j,i,k)\mid  i\in \Lambda^{(j)}\cap \Lambda^{(k)}\right\}$  is a $K$-basis of $P(j)_k$  for all $j,k\in Q_0$  (see \cite[Theorem 3.2]{P}).
\item[(2)] The module generated by  $p(j,i,i)$ is the  (unique)    submodule of $P(j)$ isomorphic to $P(i)$ for any $i\in \Lambda^{(j)}$. The map
	$f_{(j,i,i)}:P(i)\to P(j)$ is an inclusion, thus  $\left\langle p(j,i,k)\right)=\im \left(P(k)\stackrel{f_{(i,i,k)}}{\rightarrow}P(i) \stackrel{f_{(j,i,i)}}{\hookrightarrow}P(j)\right)$  is  a submodule of $ P(i)\subseteq P(j)$ (see   \cite[Remark 3.1]{P}). 
\item[(3)] We have $\im \left(P(k) \stackrel{f_{(i',i',k)}}{\rightarrow} P(i') \stackrel{f_{(i,i',i')}}{\hookrightarrow}P(i)\stackrel{f_{(j,i,i)}}{\hookrightarrow} P(j)\right)  \subseteq P(i')\subseteq P(i) \subseteq P(j)$ for  $i',i\in \Lambda^{(j)}\cap \Lambda^{(k)}$ with $i\ma i'$,  thus   $p(j,i',k)\in  P(i)_k\subseteq P(j)_k$. 
We obtain that the set  $\left\{p(j,i',k)\mid i'\in \Lambda^{(i)}\right\}$  is a $K$-basis of the space $P(i)_k$
 of  the submodule $P(i)$ of $P(j)$, since $\dimm_KP(i)_k=[P(i):S(k)]=\left|\Lambda^{(i)}\cap \Lambda^{(k)}\right|\stackrel{k\ma i}{=}\left|\Lambda^{(i)}\right|$ (see   \cite[Lemma 2.1]{P}).
\end{itemize}

We  now  consider  the  $A$-modules generated by $p(j,i,k)$ for all $i,j,k\in Q_0$ with $i\in \Lambda^{(j)}\cap \Lambda^{(k)}$. We show that   the submodule $\left\langle p(j,i,k)\right)$ of $P(j)$ is $F$-invariant for any $F\in \End_A(P(j))$, i.e.   $F\left(\left\langle p(j,i,k)\right)\right)\subseteq \left\langle p(j,i,k)\right)$. Thus   $\left\langle p(j,i,k)\right)$ is an $\End_A(P(j))^{op}$-module.

\begin{num}\begin{normalfont}\textbf{Theorem.}\end{normalfont} Let  $A=(KQ/\II,\ma)$ be a 1-quasi-hereditary algebra and  $j,k\in Q_0$.  For any  $i\in \Lambda^{(j)}\cap \Lambda^{(k)}$  the module generated by $p(j,i,k)$ is $F$-invariant for every $F\in \End_A(P(j))$. \end{num}

The proof of this Theorem follows from some properties of the lattice of  submodules of $P(j)$ generated by $p(j,i,k)$. 
 They  provide  a
  relationship    between 
 the   Hasse diagrams for the posets     $\left(\left\{\left\langle p(j,i,k) \mid i\in  \Lambda^{(j)}\cap \Lambda^{(k)}\right)\right\}, \subseteq\right)$  and $\left(\Lambda^{(j)}\cap \Lambda^{(k)}, \geqslant \right)$    as well as  $\left(\left\{\left\langle p(j,i,k) \mid k\in  \Lambda_{(i)} \right)\right\}, \subseteq\right)$  and $\left(\Lambda_{(i)}, \ma \right)$.

\begin{num}\begin{normalfont}\textbf{Lemma.}\end{normalfont} Let  $A=(KQ/\II,\ma)$ be a 1-quasi-hereditary algebra and  $j\in Q_0$.  For all  $i\in \Lambda^{(j)}$  and  $k\in \Lambda_{(i)}$ the following holds: 
\begin{itemize}
	\item[(a)]  $\left\langle p(j,i',k)\right) \subset \left\langle p(j,i,k)\right)$ if and only if $i'>i$.
	\item[(b)] $\left\langle p(j,i,k')\right) \subset \left\langle p(j,i,k)\right)$ if and only if $k'<k.$
\end{itemize} \label{ord}
\end{num}
\parbox{9.5cm}{ 
The picture to the right  
 visualizes   the presentation of $P(1)$ over the  1-quasi-hereditary algebra $A$ given by the quiver and relations in Example 1 (this algebra corresponds to a regular block of $\OO(\ger{sl}_3(\CC))$): 
 The circles   represents the  spaces $P(1)_k$  for any $k\in Q_0$  which is   spanned   by    $\left\{\textbf{p}(i,k)\mid  i\in \Lambda^{(k)}\right\}$  (see(1)), where  
\begin{center}
$\textbf{p}(i,k):=p(1,i,k)$. 
\end{center}
   The meaning of the arrows  $\textbf{p} \rightarrow \textbf{q}$ and  $\textbf{p} \dasharrow \textbf{q}$,    is   $\left\langle \textbf{p}\right) \subset \left\langle \textbf{q}\right)$ (they illustrate Lemma~\ref{ord} \textit{(a)} and \textit{(b)} respectively).  
   The set $\left\{\textbf{p}(i,k)\mid i\in \Lambda^{(j)}, \ k\in \Lambda_{(i)}\right\}$  is a $K$-basis of the submodule  $P(j)$ of $P(1)$ which is generated by $\textbf{p}(j,j)$.  The set $\left\{\textbf{p}(6,k) \mid k\in \Lambda_{(j)}\right\}$ is a $K$-basis of the  submodule $\De(j)$ of $P(1)$ which is generated by  $\textbf{p}(6,j)$ (see \cite[Remark 3.1]{P}).
}
\psset{xunit=0.52mm,yunit=0.55mm,runit=6mm}
\begin{pspicture}(0,2)(0,0)
\rput(74,-35){
\begin{tiny}
\psellipse[linewidth=1pt, linecolor=gray!90](0,80)(11,12)
\psellipse[linewidth=1pt, linecolor=gray!90](34,55)(14,15)
\psellipse[linewidth=1pt, linecolor=gray!90](-34,55)(14,15)
\psellipse[linewidth=1pt, linecolor=gray!90](34,16)(16.5,16.9)
\psellipse[linewidth=1pt, linecolor=gray!90](-34,16)(16.5,16.9)
\psellipse[linewidth=1pt, linecolor=gray!90](0,-14.5)(19,21)
\begin{scriptsize}\rput(0,108){\rnode{rasyti}{Presentation of $P(1)$}}
\rput(0,101){\rnode{rasyti2}{with the basis $\left\{p(1,i,k)\mid i,k\in Q_0, \ i\in \Lambda^{(k)}\right\}$}}\end{scriptsize}
\rput(0,80){\rnode{66o}{}}
\rput(-5,77){\rnode{66v}{}}
\rput(0,80){\rnode{666}{\textbf{p(6,6)}}}
\rput(23,88){\rnode{P6}{$P(1)_6$}}
\rput(16,88.5){\rnode{P6p}{}}
\rput(9,84){\rnode{P6g}{}}
\ncarc[arcangle=-35,linecolor=black!80]{->}{P6p}{P6g}
 
\rput(30,60){
\rput(4,0){
\rput(23,8){\rnode{P6}{$P(1)_5$}}
\rput(16,8.5){\rnode{P5p}{}}
\rput(9,4){\rnode{P5g}{}}
\ncarc[arcangle=-35,linecolor=black!80]{->}{P5p}{P5g}}
             \rput(0,0){\rnode{65o}{}}
             \rput(10,-10){\rnode{55o}{}}
             \rput(0,0){\rnode{65}{\textbf{p(6,5)}}}
             \rput(10,-10){\rnode{55}{\textbf{p(5,5)}}}
             
}
\rput(-30,60){
\rput(-4,0){
\rput(-23,8){\rnode{P4}{$P(1)_4$}}
\rput(-16,8.5){\rnode{P4p}{}}
\rput(-9,4){\rnode{P4g}{}}
\ncarc[arcangle=35,linecolor=black!80]{->}{P4p}{P4g}}
              \rput(0,0){\rnode{64o}{}}
              \rput(-10,-10){\rnode{44o}{}}
              \rput(0,0){\rnode{64}{\textbf{p(6,4)}}}
              \rput(-10,-10){\rnode{44}{\textbf{p(4,4)}}}
 
}
\rput(-20,17){
\rput(-18,0){
\rput(-23,8){\rnode{P2}{$P(1)_2$}}
\rput(-16,8.5){\rnode{P2p}{}}
\rput(-9,4){\rnode{P2g}{}}
\ncarc[arcangle=35,linecolor=black!80]{->}{P2p}{P2g}}
              \rput(-10,12){\rnode{62o}{}}
              \rput(-20,0){\rnode{42o}{}}
              \rput(-4,-2){\rnode{52o}{}}
              \rput(-4,0){\rnode{52u}{}}
              \rput(-20,-10){\rnode{22o}{}}  
              
              \rput(-10,10){\rnode{62}{\textbf{p(6,2)}}}
              \rput(-20,0){\rnode{42}{\textbf{p(4,2)}}}
              \rput(-4,-2){\rnode{52}{\textbf{p(5,2)}}}
              \rput(-20,-10){\rnode{22}{\textbf{p(2,2)}}}
               }
\rput(20,17){
\rput(18,0){
\rput(23,8){\rnode{P3}{$P(1)_3$}}
\rput(16,8.5){\rnode{P3p}{}}
\rput(9,4){\rnode{P3g}{}}
\ncarc[arcangle=-35,linecolor=black!80]{->}{P3p}{P3g}}
              \rput(10,10){\rnode{63o}{}}
              \rput(10,12){\rnode{63u}{}}
              \rput(20,0){\rnode{53o}{}}
              \rput(4,-2){\rnode{43o}{}}
              \rput(4,0){\rnode{43u}{}}
              \rput(20,-10){\rnode{33o}{}}  
              
              \rput(10,10){\rnode{63}{\textbf{p(6,3)}}}
              \rput(20,0){\rnode{53}{\textbf{p(5,3)}}}
              \rput(4,-2){\rnode{43}{\textbf{p(4,3)}}}
              \rput(20,-10){\rnode{33}{\textbf{p(3,3)}}}
        }
       
       \rput(0,2){\rnode{6l}{}} 
       \rput(0,2){\rnode{6r}{}} 
        \rput(-10,-10){\rnode{4l}{}} 
        \rput(-9,-13){\rnode{4r}{}} 
        \rput(10,-10){\rnode{5l}{}} 
        \rput(9,-13){\rnode{5r}{}} 
        \rput(-10,-20){\rnode{2l}{}} 
        \rput(10,-20){\rnode{3l}{}} 
        \rput(0,-31){\rnode{1l}{}}  
        
        \rput(2,-40){
\rput(23,8){\rnode{P2}{$P(1)_1$}}
\rput(16,9){\rnode{P1p}{}}
\rput(9,10){\rnode{P1g}{}}
\ncarc[arcangle=35,linecolor=black!80]{->}{P1p}{P1g}}
       \rput(0,2){\rnode{61}{\textbf{p(6,1)}}} 
        \rput(-10,-10){\rnode{41}{\textbf{p(4,1)}}} 
        \rput(10,-10){\rnode{51}{\textbf{p(5,1)}}}  
       \rput(-10,-20){\rnode{21}{\textbf{p(2,1)}}} 
        \rput(10,-20){\rnode{31}{\textbf{p(3,1)}}} 
        \rput(0,-30){\rnode{11}{\textbf{p(1,1)}}}   
                        
              \end{tiny}}
              \psset{nodesep=0.5pt,offset=0pt,linecolor=black,arrows=<-}
              \ncline{44}{64}
              \ncline{55}{65}
              
              \ncline{52}{62}
              \ncline{42}{62}
              \ncline{22}{52}
              \ncline{22}{42}
              \ncline{33}{53}
              \ncline{33}{43}
              \ncline{43}{63}
              \ncline{53}{63}
              \ncline{41}{61}
              \ncline[linecolor=black]{51}{61}
              \ncline{21}{51}
              \ncline{21}{41}
              \ncline{31}{51}
              \ncline{31}{41}
              \ncline{11}{21}
              \ncline{11}{31}
              \psset{nodesep=0.5pt,offset=0pt,linecolor=black!80,arrows=->}
              \ncline{pi1}{pi12}
              
              \psset{nodesep=7pt,linecolor=black!80,linestyle=dashed,arrows=->}
              \ncline{6l}{62o}
              \ncline{5r}{52u}
              \ncline{4l}{42o}
              \ncline{2l}{22o}
                \ncline{6r}{63u}
              \ncline{5l}{53o}
              \ncline{4r}{43u}
              
              \ncline{3l}{33o}
              \ncline{62o}{64o}
              \ncline{42o}{44o}
              \ncline{62o}{65o}
              \ncline{52o}{55o}

              \ncline{63o}{65o}
              \ncline{53o}{55o}
              \ncline{63o}{64o}
              \ncline{43o}{44o}
              
              \ncline{64o}{66o}
              \ncline{65o}{66o}
              \ncarc[arcangle=-35]{5}{57}
              \ncarc[arcangle=35]{v5}{v7}
              \psset{nodesep=0.5pt,offset=0pt,linecolor=blue,arrows=->}
              \ncline{pi21}{pi22}
              
\end{pspicture}
\\

\begin{num}\begin{normalfont}\textbf{Remark.} 
General theory of modules over basic algebras says that  an  $\A$-module $M$ has finitely many local submodule with $\topp$ isomorphic to $S(k)$ if and only if  there exist a $K$-basis of $\left\{x_1,x_2, \ldots, x_m\right\}$ of $M_k$ with $\left\langle x_1\right) \subset \left\langle x_2\right) \subset\cdots \subset \left\langle x_m\right)$. If $M_k$ has  such a $K$-basis  for  any $k\in Q_0(\A)$ then the number of submodules of $M$ is finite.

Lemma ~\ref{ord} implies that a projective indecomposable $A$-module $P(j)$ over 1-quasi-hereditary algebra $A$ has finitely many  submodules   if and only if the set   $(\Lambda^{(j)},  \ma )$ totally ordered. Let $A$ be an algebra given in Example 4, then 
 the  $A$-module $P(i)$ has finitely many submodules  for every  $i\in Q_0\backslash \left\{1\right\}$.  
\end{normalfont} \end{num}

The proof of  Lemma~\ref{ord} is based  on  some properties of  the  factor algebras $A(l):= A/A\left(\sum_{j\in Q_0\backslash \Lambda_{(l)}}e_j\right)A$ of $A$  considered in \cite[Sect. 5]{P}: Any $A(l)$-module is an $A$-module $M$ \\with $[M:S(t)]=0$ for all $t\in Q_0\backslash \Lambda_{(l)}$.  
The algebra $A(l)$ is quasi-hereditary   with $\De_{(l)}(k)\cong \De(k)$ for all $k\in \Lambda_{(l)}$ (the quiver of $A(l)$ is the full subquiver of $Q$ containing the vertices of  $\Lambda_{(l)}$).  Let  $p$ be a path    of $A$ passing  through the vertices in $\Lambda_{(l)}$, then this path is also a path of $A(l)$ which   we denote  by $p_{(l)}$.  The $A(l)$-module  generated by  $ p_{(l)} $ [also denoted by $\left\langle p_{(l)}\right)$]   is the largest     factor module $M$ of the  $A$-module  $\left\langle p\right)$ with the property $[M:S(t)]=0$ for all $t\in Q_0\backslash \Lambda_{(l)}$. In particular,   $\left\langle p_{(l)}\right) \subset \left\langle q_{(l)}\right)$ implies  $\left\langle p\right)\subset \left\langle q\right)$. For any $i\in \Lambda_{(l)}$ a path of the form $p(j,i,k)$ runs over the vertices  of   $\Lambda_{(l)}$ and the set $\left\{p_{(l)}(j,i,k)\mid j,i,k\in \Lambda_{(l)}, \ j,k\in \Lambda_{(i)} \right\}$ is a $K$-basis of $A(l)$ (see \cite[Lemma 5.2(a)]{P}).
\\

\textit{Proof of Lemma ~\ref{ord}.}  \textit{(a)} \dq $\Rightarrow$\dq\  
 The property (2) shows that   $ p(j,i',k)\in P(i)_k$  can be considered as an element of  a  $K$-basis  $\left\{p(j,l,k)\mid l\in \Lambda^{(i)}\right\}$ of $P(i)_k$ (we consider $P(i)$ as a submodule of $P(j)$). Thus   $i'\in \Lambda^{(i)} $ and $i'\neq i$, since   $\left\langle p(j,i',k)\right) \neq  \left\langle p(j,i,k)\right)$.

\dq $\Leftarrow$\dq\ Let $i\triangleleft l \triangleleft \cdots \triangleleft i'$. We consider the algebra $A(l)$ and show  $\left\langle p_{(l)}(j,l,k)\right)\subset \left\langle p_{(l)}(j,i,k)\right)$. This implies $\left\langle p(j,l,k)\right)\subset \left\langle p(j,i,k)\right)$ and  iteratively we have  $\left\langle p(j,i',k)\right)\subset \cdots \subset \left\langle p(j,l,k)\right) \subset \left\langle p(j,i,k)\right)$.
  
Since  $l $ is maximal in  $ \Lambda_{(l)}$, we obtain $P_{(l)}(l) = \De_{(l)}(l)\cong \De(l)$.  The equality    $[P_{(l)}(v):S(l)]=(P_{(l)}(v):\De_{(l)}(l))=1$  implies  that   $\De_{(l)}(l)$ is an unique submodule of $P_{(l)}(v)$  with top isomorphic to $S(l)$ for any $v\in \Lambda_{(l)}$  (see \cite[Lemma 5.2]{P}). Moreover,  it is a fact  that  a local submodule $M$ of $\De_{(l)}(l)$ with $\topp (M)\cong S(t)$ is unique  and  isomorphic to  $ \De_{(l)}(t)$   $(\ast)$  (see \cite[Lemma 2.5]{P}). 
 For the $A(l)$-map $f_{(j,l,k)}$ corresponding to $p_{(l)}(j,l,k)$ we have $f_{(j,l,k)}:P_{(l)}(k)\stackrel{f_{(l,l,k)}}{\longrightarrow} \De_{(l)}(l) \stackrel{f_{(j,l,l)}}{\longrightarrow}P_{(l)}(j)$  and 
  $f_{(j,l,l)} \neq 0$, since $p_{(l)}(j,l,k) \neq 0$.   Thus  $f_{(j,l,l)}$ is injective and  $\im (f_{(j,l,k)})= \left\langle p_{(l)}(j,l,k)\right)$  is a submodule of $\De_{(l)}(l)\subseteq P_{(l)}(j)$ with $\topp$ isomorphic to $S(k)$. 
 We obtain  $\left\langle p_{(l)}(j,l,k)\right) \cong \De_{(l)}(k)$ and  $\left\langle p_{(l)}(j,i,k)\right) \subseteq P_{(l)}(i) \backslash \De_{(l)}(l)$ because  $p_{(l)}(j,i,k)$ and $p_{(l)}(j,l,k)$ are linearly independent.
 
   Let $W:=\De_{(l)}(l)\cap \left\langle p_{(l)}(j,i,k)\right)$.
  We show now existence of  $t\in \Lambda_{(l)}$   with $t\triangleright k$ such that $S(t)$ is a direct summand of $\topp \left(W\right)$. This implies an existence of a local submodule $L$ of $\De_{(l)}(l)$ with $\topp L\cong S(t)$ such that $L\subseteq \left\langle p_{(l)}(j,i,k)\right)$.  The fact $(\ast)$ provides  $L=\De_{(l)}(t)$ and  consequently   $\left\langle p_{(l)}(j,l,k)\right) = \De_{(l)}(k) \subset \De_{(l)}(t) \subseteq \left\langle p_{(l)}(j,i,k)\right)$ (see \cite[Remark 2.6]{P}):  
The filtration $0\subset \De_{(l)}(l)\subset P_{(l)}(i)$ is $\De$-good, since  $\left(P_{(l)}(i):\De_{(l)}(j)\right)=1$ for $j=i,l$ (see \cite[Lemma 5.2]{P}). Since  $\left\langle p_{(l)}(j,i,k)\right)/W  \hookrightarrow P_{(l)}(i)/\De_{(l)}(l) \cong \De_{(l)}(i)$, we have $\left\langle p_{(l)}(j,i,k)\right)/W \cong \De_{(l)}(k)$ (see $(\ast)$).  For  the  $A$-module   $\left\langle p_{(l)}(j,i,k)\right)$   there exists a surjective map $F:P(k) \twoheadrightarrow \left\langle p_{(l)}(j,i,k)\right)$. According to \cite[Lemma 2.4]{P} we have  $P(k)/\left(\sum_{k\triangleleft t}P(t)\right) \cong \De(k)$, thus  $F\left(\sum_{k\triangleleft t}P(t)\right) = W$.   
 Since   $W$ is an $A(l)$-module   we obtain     
     $\topp \left(W\right)\in  \add \left(\topp \left(\sum_{k\triangleleft t}P(t)\right)\right)\cap \add \left(\bigoplus_{t\in \Lambda_{(l)}}S(t)\right) = \add \left(\bigoplus_{t\in \Lambda_{(l)}\atop k\triangleleft t}S(t)\right)$ (see  \cite[Corollary 3.9]{ASS}).

 \textit{(b)} $\dq \Leftrightarrow\dq$ It is  $\left\langle p(j,l,k')\right) \subset \left\langle p(j,l,k)\right)$ if and only if  $\left\langle  p_{(l)}(j,l,k')\right) \subset \left\langle p_{(l)}(j,l,k)\right)$.  We have    $\left\langle p_{(l)}(j,l,k')\right)\cong \De_{(l)}(k')\cong \De(k')$ and  $\left\langle  p_{(l)}(j,l,k)\right)\cong \De_{(l)}(k)\cong \De(k)$. 
 According to  \cite[Lemma 2.5]{P} we obtain  $\De (k')\subset \De(k)$ if and only if $k'<k$. \hfill $\Box$
\\

\textit{Proof of Theorem.} 
Let $F\in \End_A(P(j))$.  The submodule $P(i)$ of $P(j)$ is $F$-invariant according to \cite[Remark 2.3]{P}.  
 Since $\left\langle p(j,i,k)\right) \subseteq P(i)\subseteq P(j)$, we have   $F\left( p(j,i,k)\right) \in P(i)_k  \stackrel{\text{(3)}}{=}\spann_K\left\{p(j,i',k)\mid i'\in \Lambda^{(i)}\right\}$.   Lemma ~\ref{ord} \textit{(a)} implies $F(p(j,i,k))\in \left\langle p(j,i,k)\right)$, thus $F(\left\langle p(j,i,k)\right)) \subseteq \left\langle p(j,i,k)\right)$. 
\hfill $\Box$

\section{Generators of the characteristic tilting module}
In this section  we consider a 1-quasi-hereditary algebra $A=KQ/\II$  for  which  the corresponding Ringel-dual $R(A)$  is also 1-quasi-hereditary. According to \cite[Theorem 6.1]{P} this is  exactly the case if 
the factor algebra  $A(i):=A/\left[A\left(\sum_{t\in Q_0\backslash \Lambda_{(i)}}e_t\right)A\right]$ is  1-quasi-hereditary for every $i\in Q_0$ (for $A(i)$-modules we use  the index $(i)$). 

The  properties of a 1-quasi-hereditary algebra $(\A,\ma )$ (with $\left\{1\right\}=\min\{Q_0(\A),\ma\}$) yields  the following: If $x\in \A$ generates $P_{\A}(1)$, then  $p(1,k,k)\cdot x$  generates the submodule $P_{\A}(k)$ of $P_{\A}(1)$ for any increasing path   $p(1,k,k)$   from $1$ to $k$ (see \cite[Remark 3.1]{P}). 

Let  $A(i)$ be  1-quasi-hereditary for all $i\in Q_0$. The direct summand  $T(i)$ of the characteristic tilting $A$-module $T$  is a local submodule of $P(1)$ with $\topp T(i)\cong S(1)$ (see \cite[Theorem 5.1]{P}).  
There exists  $\ger{t}(i, 1)\in P(1)_1$ which  generates $T(i)$ for all $i\in Q_0$.     Because 
 $i\in Q_0(A(i))=\Lambda_{(i)}$ is maximal,  for $A(i)$-module $T(i)$  we have   
 $P_{(i)}(1)\cong I_{(i)}(1) \cong T(i)\cong T_{(i)}(i)$  (see \cite[Lemma 5.2]{P}). 
 Since $P_{(i)}(1)\cong   \left\langle \ger{t}(i, 1)\right)$, for any  $k\in \Lambda_{(i)}$ we obtain   that 
\begin{center}
$\ger{t}(i, k):=p(1,k,k)\cdot \ger{t}(i, 1)$ \   generates  \   the submodule  \   $P_{(i)}(k)$  \   of $P_{(i)}(1)$.
\end{center}
Since $\ger{t}(i, k)\in P(1)_k$,  we have   $\ger{t}(i, k)= \sum_{l\in \Lambda^{(k)}}c_l\cdot p(1,l,k)$  for some $c_l\in K$ (see Section 2 (1)).
In particular,
 for any $j\in \Lambda^{(i)}$ the algebra $A(i)$ is  a  factor algebra of $A(j)$, namely   $A(i)=\big(A(j)\big)(i)$, thus 
  $T(i)$  is a direct summand of the characteristic tilting $A(j)$-module $T_{(j)}=\bigoplus_{k\in \Lambda_{(j)}}T_{(j)}(k)$  with $T(i)\cong T_{(j)}(i)$.
The next Theorem shows particular 	features of $\ger{t}(i,k)$ and of $A$-modules  generated by $\ger{t}(i,k)$.

\begin{num}\begin{normalfont}\textbf{Theorem.}\end{normalfont} Let  $A$ and $R(A)$  be 1-quasi-hereditary and     $\ger{t}(i,k)$   be an element    defined as above.   Then the following hold:  
\begin{itemize}
	\item[(1)] The set  $\left\{ \ger{t}(i,k)\mid i\in Q_0, \  k \in \Lambda_{(i)} \right\} $   is  a   $K$-basis of  $P(1)$.
	\item[(2)] $ \left\{\left\langle \ger{t}(i,k)\right)\mid  i\in Q_0, \  k\in \Lambda_{(i)}\right\} = \left\{\text{local submodules  of   } P(1)  \text{  in   } \ger{F}(\De) \right\}.$
	\\[3pt]  Let  $j\in Q_0 $, then  a local submodule of $P(j)$ in  $\ger{F}(\De)$  is isomorphic to some   $ \left\langle \ger{t}(i,k)\right)$. 
	\item[(3)] Let $F\in \End_A(P(1))$, then  $\left\langle \ger{t}(i,k)\right)$ is $F$-invariant  for  all $i,k\in Q_0$ with $k\in \Lambda_{(i)}$.
\end{itemize}
 \label{basistwo} \end{num}

\textit{Proof.} 
 \textit{(1)} Obviously, $\ger{t}(i, k)$   belongs to $P(1)_k$ for all  $i\in \Lambda^{(k)}$.
  Let $r = \left|\Lambda^{(k)}\right|$  and   $ \mathcal{L}(k):=\left\{(i_1, i_2, \ldots , i_r) \mid i_l\in \Lambda^{(k)} , \   i_{v}\not\geqslant i_{t}, \ 1\leq  v < t \leq r\right\}$   and $(i_1, \ldots , i_t,  \ldots , i_r)\in \mathcal{L}(k)$.   By induction on $t$ we  show  that  $\ger{t}(i_1,k), \ldots , \ger{t}(i_r,k)$ are  linearly  independent. This    implies that the set  $\left\{\ger{t}(i,k) \mid i\in \Lambda^{(k)} \right\}$ is a $K$-basis of $P(1)_k$, since $\dimm_KP(i)_k= r$ (see \cite[Lemma 2.1]{P}):   
   For  $t=1$ we have $i_1=k$, thus    $\ger{t}(i_1,i_1)\neq 0$ is linearly  independent.   By  definition of $\mathcal{L}(k)$    we have $i_l\not\geqslant i_{t+1}$,  thus $ [\left\langle \ger{t}(i_{l},k)\right):S(i_{t+1})]=[P_{(i_{l})}(k):S(i_{t+1})]=0$  for all  $1\leq l\leq t$,  however  $[\left\langle \ger{t}(i_{t+1},k)\right):S(i_{t+1})]=[P_{(i_{t+1})}(k) :S(i_{t+1})]=1$, since $A(i_{t+1})$ is 1-quasi-hereditary.  Hence $\left\langle \ger{t}(i_{t+1},k)\right)\not\subseteq \sum_{l=1}^{t} \left\langle  \ger{t}(i_l,k)\right)$ and consequently $\ger{t}(i_{t+1},k) \not\in \spann_K\left\{ \ger{t}(i_1,k), \ldots , \ger{t}(i_{t},k)\right\}$.

\textit{(2)} $\dq \supseteq\dq$ 
 Let    $M\in \ger{F}(\De)$ with $M\subseteq P(1)$  and $\topp M\cong S(k)$. Then 
 $M\cong P(k)/M'$ with $M'=D(t+1)=\sum_{m=t+1}^{r}P(i_m)$ for some $\textbf{i}=(i_1,\ldots, i_t,i_{t+1},\ldots, i_r)\in \mathcal{L}(k)$ (see \cite[Proposition 4.2]{P}).  The $\De$-good filtration  $\mathscr{D}(\textbf{i})$  of $P(k)$ induces  the $\De$-good filtration $0\subset D(t)/M'\subset \cdots \subset D(2)/M'\subset D(1)/M'=P(k)/M' = M$  with $D(t)/M'\cong \De(i_{t})$ and $M/(D(2)/M')\cong \De(k)=\De(i_1)$.  Since  $\soc M$ and $\topp M$ are  simple, any $\De$-good filtration  starts with $\De(i_{t})$ and the top  quotient   is   $\De(i_1)$, moreover $\left(M:\De(l)\right)=\begin{small}\left\{
\begin{array}{cl}
1 & \text{if   } l\in \left\{i_1,\ldots , i_t\right\},\\
0 & \text{if } l\in \left\{i_{t+1},\ldots , i_r\right\}.
\end{array}
\right. \end{small}$   Definition of $\mathcal{L}(k)$  implies    $i_1\ma l\ma i_{t}$  for any $l\in \left\{i_1,\ldots , i_{t}\right\}$ and $\Lambda_{(i_t)}\cap \left\{i_{t+1},\ldots , i_r\right\}=\emptyset$. Thus  $\left\{i_1,\ldots , i_{t}\right\}=\Lambda^{(k)}\cap \Lambda_{(i_{t})}$ and   $\left\{i_{t+1},\ldots , i_r\right\}=\Lambda^{(k)}\backslash  \Lambda_{(i_t)}$.  According to   \cite[Lemma 5.2]{P}  we have   $M\cong P(k)/\left(\sum_{l\in \Lambda^{(k)}\backslash \Lambda_{(i_t)}}P(l)\right)\cong P_{(i_t)}(k)$, thus $M= \left\langle \ger{t}(i_t,k)\right) $. 
 Since    $ \ger{F}(\De_{(i_t)}) \subseteq \ger{F}(\De)$,  we have $\dq \subseteq\dq$. 
 
Any  submodule of $P(j)$  is isomorphic to a submodule of  $P(1)$, since $P(j)\hookrightarrow P(1)$ for all $j\in Q_0$. Thus any local submodule of $P(j)$ in  $\ger{F}(\De)$ is  a submodule of  $P(1)$ in $\ger{F}(\De)$.
 
 \textit{(3)}  First, we show that any $A(i)$-submodule $M$ of $P(1)$ is contained in $\left\langle \ger{t}(i,1) \right) \cong P_{(i)}(1)\cong T(i)$: Let $i\neq n$. We have  $P(1)/T(i)\in \ger{F}(\nabla)$, since the category $\ger{F}(\nabla)$ is closed under  cokernels of injective maps (see \cite{Rin1}).  If   $(P(1)/T(i):\nabla(t))\neq 0$,  then $t\in Q_0\backslash \Lambda_{(i)}$ because 
   $(P(1):\nabla(t))=1$ for  all $t\in Q_0$ and   $(T(i):\nabla(t))=(I_{(i)}(1):\nabla(t))=1$ for all $t\in \Lambda_{(i)}$.  Thus $\soc \left(P(1)/T(i)\right)\in \add \left\{\bigoplus_{t\in  Q_0\backslash \Lambda_{(i)}}S(t)\right\}$. Consequently, for any submodule $N$ of $P(1)$ with $ \left\langle \ger{t}(i,1) \right) \cong T(i) \subset N $ we have  $[N:S(t)]\neq 0$ for some $t\in Q_0\backslash \Lambda_{(i)}$, i.e. $N$ is not a $A(i)$-module. Thus  $ \left\langle \ger{t}(i,1) \right)$ is the largest $A(i)$-submodule of $P(1)$.
 
 Let $M$ be a submodule of $P(1)$ with  $[M:S(t)]=0$  for all $t\in Q_0\backslash \Lambda_{(i)}$, then for any $F\in \End_A(P(1))$  it is     $[F(M):S(t)]=0$ for  $t\in Q_0\backslash \Lambda_{(i)}$. 
 We have $F(M) \subseteq \left\langle \ger{t}(i,1)\right) $. Since \\  $A(i)$ is 1-quasi-hereditary,  any submodule $L$ of $P_{(i)}(1)$ with $\topp L\cong S(k)$ is contained  in  $ P_{(i)}(k) \cong \left\langle \ger{t}(i,k)\right)$ (see \cite[Lemma 2.2]{P}). Thus  $F(\left\langle \ger{t}(i,k)\right)) \subseteq \left\langle \ger{t}(i,k)\right)$ for every $F\in \End_A(P(1))$, since  $F(\left\langle \ger{t}(i,k)\right) ) $ is a submodule of $P_{(i)}(1)$ with  $\topp F(\left\langle \ger{t}(i,k)\right) ) \cong S(k)$.
  \hfill $\Box$
\\

A similar relationship to  that between the modules generated by $p(1,i,k)$ and the posets $(\Lambda_{(i)},\ma)$ as well as $(\Lambda^{(k)},\geqslant)$  in Lemma~\ref{ord}  exists  also between  the submodules  generated by $\ger{t}(i, k)$ and $(\Lambda^{(k)},\ma)$ as well as $(\Lambda_{(i)},\geqslant)$: The Hasse diagrams of  the    sets $\left(\left\{\left\langle \ger{t}(i, k)\right)\mid i \in \Lambda^{(k)}\right\},\subseteq \right)$   and $\left( \Lambda^{(k)}, \ma \right)$  as well as the diagrams  of  $\left(\left\{\left\langle \ger{t}(i, k)\right)\mid k\in \Lambda_{(i)}\right\},\subseteq\right)$  and $\left(\Lambda_{(i)}, \geqslant \right)$   coincide.

\begin{num}\begin{normalfont}\textbf{Lemma.}\end{normalfont} Let  $A$ and $R(A)$  be 1-quasi-hereditary and     $\ger{t}(i,k)$   be an element    defined as above. Then  for all $i,i',k,k'\in Q_0$ it is 
\begin{itemize}
	\item[(a)] $\left\langle \ger{t}(i',k)\right) \subset \left\langle \ger{t}(i,k)\right)$ if and only if $i'<i$.
	\item[(b)] $\left\langle \ger{t}(i,k')\right) \subset \left\langle \ger{t}(i,k)\right)$ if and only if $k'>k$.
\end{itemize}
\label{letzte}
 \end{num}

\textit{Proof.} \textit{(a)}  If  $\underbrace{\left\langle \ger{t}(i',k)\right)}_{P_{(i')}(k)} \subset \underbrace{\left\langle \ger{t}(i,k)\right)}_{ P_{(i)}(k)}$,     then  $\left[P_{(i')}(k):S(i')\right]\neq 0$. Thus $\left[P_{(i)}(k):S(i')\right]\neq 0$ so we get $i'\in \Lambda_{(i)}$ (see \cite[Lemma 5.2 (c)]{P}).    Since   $\left\langle \ger{t}(i',k)\right) \neq  \left\langle \ger{t}(i,k)\right)$, we obtain  $i\neq i'$.
On the other hand,  if  $i'< i$, then  $T(i')$ is the  direct summand $T_{(i)}(i')$ of the characteristic tilting $A(i)$-module  and therefore a submodule of $P_{(i)}(1)$, since $A(i')=\big(A(i)\big)(i')$ is 1-quasi-hereditary  (see \cite[Theorem 5.1]{P}).  We have
 $P_{(i')}(k) \hookrightarrow P_{(i')}(1) \cong  T(i') \cong T_{(i)}(i')\hookrightarrow  P_{(i)}(1)$, thus $P_{(i')}(k)$ is a local submodule of $ P_{(i)}(1)$  with top isomorphic to  $S(k)$. According  to \cite[Lemma 2.2]{P},  we obtain  $P_{(i')}(k) \subseteq  P_{(i)}(k)$.
 Since $i'<i$ we have  $\left[P_{(i')}(k):S(i)\right]=0$ and $\left[P_{(i)}(k): S(i)\right]=1$. Therefore  $P_{(i')}(k) \neq P_{(i)}(k)$  and consequently  $\left\langle \ger{t}(i',k)\right) \subset \left\langle \ger{t}(i,k)\right)$.

\textit{(b)} According to \cite[Lemma 2.2]{P}  we have  $P_{(i)}(k') \subset P_{(i)}(k)$ if and only if $k'>k$. \hfill $\Box$
\\

  Algebras  associated to  the blocks of category $\OO(\ger{g})$    are  Ringel self-dual (see ~\cite{So}). The 1-quasi-hereditary algebra $A$ given by quiver and relations in the Example 1 corresponds to the regular block of $\OO(\ger{sl}_3(\CC))$, thus the Ringel-dual  of $A$ is also 1-quasi-hereditary. Analogously  to the previous picture we can give  the  inclusion diagram of the submodules $P_{(i)}(k)\cong \left\langle \ger{t}(i,k)\right)$  of $P(1)$, in view of    Lemma~\ref{letzte}. The  meaning of  circles and  arrows is the 
  \\[2pt]
\parbox{9.5cm}{same  as 
     in the last picture (the arrows $\ger{p}\rightarrow \ger{q}$ and $\ger{p} \dasharrow \ger{q}$ illustrates Lemma~\ref{letzte}\textit{(a)} and \textit{(b)} respectively).  We additionally   pointed out the elements  which generates the direct summands of the characteristic tilting module, projective and standard submodules of $P(1)$.  The notation $M \leftrightsquigarrow  \textbf{p} $  resp. $\textbf{p} \leftrightsquigarrow M$    means    $M=\left\langle \textbf{p}\right) $.
 These $\dimm_KP(1)=19$ submodules of $P(1)$ are the local submodules of $P(1)$  from $\ger{F}(\De)$.

The element $\ger{t}(i,k)$ is the following  linear combination of $\left\{\textbf{p}(l,k):=p(1,l,k)\mid l\in \Lambda^{(k)}\right\}$:

 $
\begin{array}{l}
\ger{t}(2,1)=\textbf{p}(5,1), \\
\ger{t}(3,1)=\textbf{p}(4,1), \\
\ger{t}(4,1)=\textbf{p}(3,1), \\
\ger{t}(5,1)=\textbf{p}(2,1), \\
\end{array}
$    \hspace{4mm}
$
\begin{array}{l}
\ger{t}(4,3)=\textbf{p}(5,3),\\
\ger{t}(5,3)=\textbf{p}(4,3)+ \textbf{p}(5,3), \\[6pt]
\ger{t}(4,2)=\textbf{p}(4,2)+ \textbf{p}(5,2), \\
\ger{t}(5,2)=\textbf{p}(4,2), \\
\end{array}
$
\\
$
\begin{array}{l}
\ger{t}(i,i)=\textbf{p}(6,i), \\
\ger{t}(6,i)=\textbf{p}(i,i)  \\
\end{array}
$ \   for every  $i\in Q_0$.
}
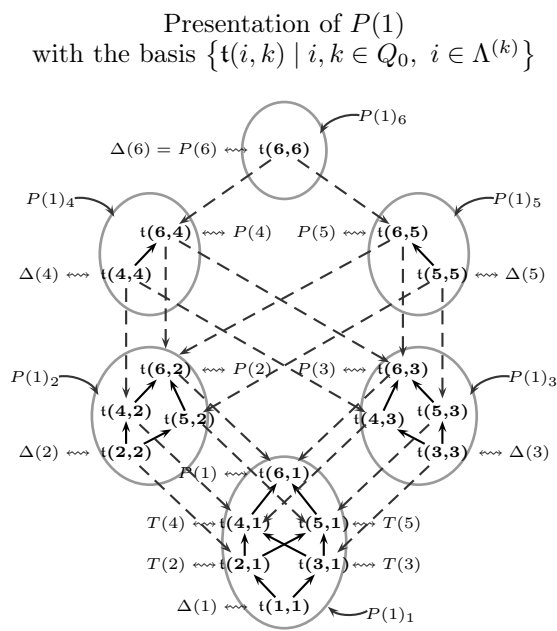
\psset{xunit=0.52mm,yunit=0.55mm,runit=6mm}
\begin{pspicture}(0,2)(0,0)
\rput(74,-39){
\begin{tiny}
\psellipse[linewidth=1pt, linecolor=gray!80](0,80)(11,12)
\psellipse[linewidth=1pt, linecolor=gray!80](34,55)(13,15)
\psellipse[linewidth=1pt, linecolor=gray!80](-34,55)(13,15)
\psellipse[linewidth=1pt, linecolor=gray!80](34,16)(15,16.9)
\psellipse[linewidth=1pt, linecolor=gray!80](-34,16)(15,16.9)
\psellipse[linewidth=1pt, linecolor=gray!80](0,-14.5)(16,21)
\begin{footnotesize}\rput(0,110){\rnode{rasyti}{Presentation of $P(1)$}}
\rput(0,103){\rnode{rasyti2}{with the basis $\left\{\ger{t}(i,k)\mid i,k\in Q_0, \ i\in \Lambda^{(k)}\right\}$}} \end{footnotesize}
\rput(0,80){\rnode{66o}{}}
\rput(-5,77){\rnode{66v}{}}
\rput(0,80){\rnode{666}{$\ger{t}$\textbf{(6,6)}}}
\rput[r](-9,80){\rnode{x3}{$\De(6)= P(6)\leftrightsquigarrow$}}
\rput(25,88){\rnode{P6}{$P(1)_6$}}
\rput(19,88.5){\rnode{P6p}{}}
\rput(9,84){\rnode{P6g}{}}
\ncarc[arcangle=-35,linecolor=black!80]{->}{P6p}{P6g}

\rput(30,60){
\rput(4,0){
\rput(25,8){\rnode{P6}{$P(1)_5$}}
\rput(19,8.5){\rnode{P5p}{}}
\rput(9,4){\rnode{P5g}{}}
\ncarc[arcangle=-35,linecolor=black!80]{->}{P5p}{P5g}}
             \rput(0,0){\rnode{65o}{}}
             \rput(10,-10){\rnode{55o}{}}
             \rput(0,0){\rnode{65}{$\ger{t}$\textbf{(6,5)}}}
             \rput[r](-9,0){\rnode{65b}{$P(5)\leftrightsquigarrow$}}
             \rput(10,-10){\rnode{55}{$\ger{t}$\textbf{(5,5)}}}
             \rput[l](18,-10){\rnode{55b}{$\leftrightsquigarrow \De(5)$}}

}
\rput(-30,60){
\rput(-4,0){
\rput(-25,8){\rnode{P4}{$P(1)_4$}}
\rput(-19,8.5){\rnode{P4p}{}}
\rput(-9,4){\rnode{P4g}{}}
\ncarc[arcangle=35,linecolor=black!80]{->}{P4p}{P4g}}
              \rput(0,0){\rnode{64o}{}}
              \rput(-10,-10){\rnode{44o}{}}
              \rput(0,0){\rnode{64}{$\ger{t}$\textbf{(6,4)}}}
              \rput[l](9,0){\rnode{64b}{$\leftrightsquigarrow P(4)$}}
              \rput(-10,-10){\rnode{44}{$\ger{t}$\textbf{(4,4)}}}
              \rput[r](-19,-10){\rnode{44b}{$ \De(4)\leftrightsquigarrow$}}

}
\rput(-20,17){
\rput(-18,0){
\rput(-25,8){\rnode{P2}{$P(1)_2$}}
\rput(-19,8.5){\rnode{P2p}{}}
\rput(-9,4){\rnode{P2g}{}}
\ncarc[arcangle=35,linecolor=black!80]{->}{P2p}{P2g}}
              \rput(-10,12){\rnode{62o}{}}
              \rput(-20,0){\rnode{42o}{}}
              \rput(-4,-2){\rnode{52o}{}}
              \rput(-4,0){\rnode{52u}{}}
              \rput(-20,-10){\rnode{22o}{}}

              \rput(-10,10){\rnode{62}{$\ger{t}$\textbf{(6,2)}}}
              \rput[l](-1,10){\rnode{62b}{$\leftrightsquigarrow P(2)$}}
              \rput(-20,0){\rnode{42}{$\ger{t}$\textbf{(4,2)}}}
              \rput(-4,-2){\rnode{52}{$\ger{t}$\textbf{(5,2)}}}
              \rput(-20,-10){\rnode{22}{$\ger{t}$\textbf{(2,2)}}}
              \rput[r](-29,-10){\rnode{22b}{$\De(2)\leftrightsquigarrow$}}
               }
\rput(20,17){
\rput(18,0){
\rput(25,8){\rnode{P3}{$P(1)_3$}}
\rput(19,8.5){\rnode{P3p}{}}
\rput(9,4){\rnode{P3g}{}}
\ncarc[arcangle=-35,linecolor=black!80]{->}{P3p}{P3g}}
              \rput(10,10){\rnode{63o}{}}
              \rput(10,12){\rnode{63u}{}}
              \rput(20,0){\rnode{53o}{}}
              \rput(4,-2){\rnode{43o}{}}
              \rput(4,0){\rnode{43u}{}}
              \rput(20,-10){\rnode{33o}{}}

              \rput(10,10){\rnode{63}{$\ger{t}$\textbf{(6,3)}}}
              \rput[r](1,10){\rnode{62b}{$P(3)\leftrightsquigarrow $}}
              \rput(20,0){\rnode{53}{$\ger{t}$\textbf{(5,3)}}}
              \rput(4,-2){\rnode{43}{$\ger{t}$\textbf{(4,3)}}}
              \rput(20,-10){\rnode{33}{$\ger{t}$\textbf{(3,3)}}}
              \rput[l](29,-10){\rnode{33b}{$\leftrightsquigarrow \De(3)$}}
        }

       \rput(0,2){\rnode{6l}{}}
       \rput(0,2){\rnode{6r}{}}
        \rput(-10,-10){\rnode{4l}{}}
        \rput(-9,-13){\rnode{4r}{}}
        \rput(10,-10){\rnode{5l}{}}
        \rput(9,-13){\rnode{5r}{}}
        \rput(-10,-20){\rnode{2l}{}}
        \rput(10,-20){\rnode{3l}{}}
        \rput(0,-31){\rnode{1l}{}}

        \rput(2,-40){
\rput(25,8){\rnode{P2}{$P(1)_1$}}
\rput(19,9){\rnode{P1p}{}}
\rput(9,10){\rnode{P1g}{}}
\ncarc[arcangle=35,linecolor=black!80]{->}{P1p}{P1g}}
       \rput(0,2){\rnode{61}{$\ger{t}$\textbf{(6,1)}}}
       \rput[r](-9,2){\rnode{61b}{$P(1)\leftrightsquigarrow $}}
        \rput(-10,-10){\rnode{41}{$\ger{t}$\textbf{(4,1)}}}
        \rput(10,-10){\rnode{51}{$\ger{t}$\textbf{(5,1)}}}
       \rput(-10,-20){\rnode{21}{$\ger{t}$\textbf{(2,1)}}}
        \rput(10,-20){\rnode{31}{$\ger{t}$\textbf{(3,1)}}}
        \rput(0,-30){\rnode{11}{$\ger{t}$\textbf{(1,1)}}}
        \rput[r](-9,-30){\rnode{11b}{$\De(1)\leftrightsquigarrow $}}
        \rput[r](-17,-20){\rnode{t2b}{$T(2)\leftrightsquigarrow $}}
        \rput[r](-17,-10){\rnode{t4b}{$T(4)\leftrightsquigarrow $}}
        \rput[l](17,-10){\rnode{t5b}{$\leftrightsquigarrow T(5)$}}
        \rput[l](17,-20){\rnode{t5b}{$\leftrightsquigarrow T(3)$}}

              \end{tiny}}
              \psset{nodesep=0.5pt,offset=0pt,linecolor=black,arrows=->}
              \ncline{44}{64}
              \ncline{55}{65}

              \ncline{52}{62}
              \ncline{42}{62}
              \ncline{22}{52}
              \ncline{22}{42}
              \ncline{33}{53}
              \ncline{33}{43}
              \ncline{43}{63}
              \ncline{53}{63}
              \ncline{41}{61}
              \ncline[linecolor=black]{51}{61}
              \ncline{21}{51}
              \ncline{21}{41}
              \ncline{31}{51}
              \ncline{31}{41}
              \ncline{11}{21}
              \ncline{11}{31}
              \psset{nodesep=0.5pt,offset=0pt,linecolor=black!80,arrows=<-}
              \ncline{pi1}{pi12}

              \psset{nodesep=5pt,linecolor=black!80,linestyle=dashed,arrows=<-}
              \ncline{6l}{62o}
              \ncline{5r}{52u}
              \ncline{4l}{42o}
              \ncline{2l}{22o}
                \ncline{6r}{63u}
              \ncline{5l}{53o}
              \ncline{4r}{43u}

              \ncline{3l}{33o}
              \ncline{62o}{64o}
              \ncline{42o}{44o}
              \ncline{62o}{65o}
              \ncline{52o}{55o}

              \ncline{63o}{65o}
              \ncline{53o}{55o}
              \ncline{63o}{64o}
              \ncline{43o}{44o}

              \ncline{64o}{66o}
              \ncline{65o}{66o}
              \ncarc[arcangle=-35]{5}{57}
              \ncarc[arcangle=35]{v5}{v7}
              \psset{nodesep=0.5pt,offset=0pt,linecolor=blue,arrows=<-}
              \ncline{pi21}{pi22}

\end{pspicture}

\text{  }
\\

\textbf{Acknowledgments.}  
I would like to thank  Julia Worch   for     comments and remarks.

\begin{small}

\end{small}

\end{document}